\documentclass[english,11pt,3p,number,sort&compress, export]{elsarticle}
\usepackage[T1]{fontenc}
\usepackage[latin9]{inputenc}
\usepackage{geometry}
\usepackage{color}
\definecolor{MyRed}{rgb}{0,0,0}

\usepackage{array}
\usepackage{float}
\usepackage{bm}
\usepackage{algorithm2e}
\SetArgSty{textnormal}
\usepackage{amsmath}
\usepackage{amssymb}
\usepackage{stmaryrd}
\usepackage{graphicx}
\usepackage{lipsum}

\providecommand{\tabularnewline}{\\}

\usepackage{natbib}
\usepackage{eucal}
\usepackage{dsfont}
\usepackage{comment}
\usepackage{ifthen}    
\usepackage{lscape}
\usepackage{graphicx}

\usepackage{amsmath}

\usepackage{stmaryrd} 

\usepackage{hyperref}
\usepackage{bm}

\usepackage[subpreambles=true]{standalone}
\usepackage{xspace}
\usepackage[percent]{overpic}

\usepackage[customcolors]{hf-tikz}
\usepackage{tikz}
\usepackage{pgfplots}
\usetikzlibrary{calc,shadings,patterns,tikzmark, plotmarks, spy, 
pgfplots.polar, external, matrix, shapes.symbols,shadings,shapes, 
decorations.shapes,decorations.pathmorphing,fit,backgrounds}
\pgfplotsset{compat=1.10}   

\usepgfplotslibrary{groupplots}
\usetikzlibrary{calc}
\usepackage{pgfplotstable}

\usepackage{xcolor}
\usepackage{adjustbox} 

\newcommand{\hdivpure}{\text{H}(\text{div})}
\newcommand{\hcurl}{\text{H}(\text{curl})}
\newcommand{\hone}{\text{H}^1}

\newenvironment{codingtime}{\ttfamily}{\par}

\pgfplotsset{compat=newest}

\hypersetup{urlcolor=blue, colorlinks=true}

\usepackage{amssymb}

\graphicspath{ {./figures/} }

\makeatletter
\def\input@path{{./figs_pgfplots/}{./figures/}}
\makeatother

\makeatletter
\@ifundefined{showcaptionsetup}{}{%
 \PassOptionsToPackage{caption=false}{subfig}}
\usepackage{subfig}
\makeatother

\usepackage{babel}

\journal{Book series "Advances in Applied Mechanics, Vol 59, by Elsevier\cite{ref-pattern}.}

\makeatletter
\def\ps@pprintTitle{%
     \let\@oddhead\@empty
     \let\@evenhead\@empty
     \def\@oddfoot{\footnotesize\itshape
      Preprint accepted for publication in the \ifx\@journal\@empty Elsevier
      \else\@journal\fi\hfill} 
     \let\@evenfoot\@oddfoot}
\makeatother

\begin{document}

\begin{frontmatter}{

\title{Extending h adaptivity with refinement patterns}

\author[uni]{Giovane Avancini}

\ead{giovanea@unicamp.br}

\author[uni]{Nathan Shauer}

\ead{shauer@unicamp.br}

\author[uni]{Francisco T. Orlandini}

\ead{francisco.orlandini@gmail.com }

\author[puc]{Paulo Cesar A. Lucci}

\ead{cesar.lucci@mac.com}

\author[uni]{Philippe R. B. Devloo}

\ead{phil@unicamp.br}

\address[uni]{Universidade Estadual de Campinas, R. Josiah Willard Gibbs 85 - Cidade Universitaria, Campinas SP, Brazil, CEP 13083-839}

\address[puc]{Exact Sciences, Environmental and Technologies Center, Pontifical Catholic University of Campinas (PUC-Campinas), Campinas 13086-099, Brazil}

\begin{abstract}

This contribution introduces the idea of refinement patterns for the generation of optimal meshes in the context of the Finite Element Method. The main idea is to generate a library of possible patterns on which elements can be refined and use this library to inform an h adaptive code on how to handle complex refinements in regions of interest. There are no restrictions on the type of elements that can be refined, and the patterns can be generated for any element type. The main advantage of this approach is that it allows for the generation of optimal meshes in a systematic way where, even if a certain pattern is not available, it can easily be included through a simple text file with nodes and sub-elements. The contribution presents a detailed methodology for incorporating refinement patterns into h adaptive Finite Element Method codes and demonstrates the effectiveness of the approach through mesh refinement of problems with complex geometries.

\end{abstract}

}\end{frontmatter}

\section{Introduction \label{sec:intro}}

The Finite Element Method (FEM) is a powerful numerical technique widely used in various fields of engineering and science for solving partial differential equations. One of the key advantages of the FEM is its ability to handle complex geometries and adaptively refine the mesh to accurately capture the solution behavior in regions of interest. In particular, adaptive h-p refinement, which combines both mesh refinement (h-refinement) and polynomial degree refinement (p-refinement), has been proven to be highly effective in achieving accurate and efficient solutions \cite{Babuska1988,Babuska1994}.

Overall, h-refinement combined with adaptive algorithms can play a crucial role in improving the accuracy of the solution by refining the mesh in regions of interest, especially in regions with strong gradients, boundary layers, and singularities \cite{Demkowicz2002}. Optimal mesh refinements for specific problems have been studied in the literature \cite{Babuska1986}.

In the context of h-refinement, the generation of hanging nodes is a common issue that arises when refining the mesh towards a region of interest. The meshes created in this way are often called 1-irregular meshes. Hanging nodes occur when a refined element shares a face with a coarser element, resulting in a discontinuity in the solution across the interface that needs to be treated. Methodologies to do so have been proposed in the literature with an increased interest during the 80s and 90s. The most common approach is to impose the continuity between elements as a restriction to the shape functions of both elements at a given interface \cite{DEMKOWICZ1985,ODEN1986,DEVLOO1988,DEVLOO1987}. This methodology has been implemented in the \texttt{NeoPZ} library \cite{Devloo1997} for 1D, 2D, and 3D elements of any polynomial order. Other methodologies have also been explored in the literature as for instance in \cite{Morton1995} where txransition elements are proposed that do not rely on constraints to produce interelement continuity. Recently, there has been renewed interest in the topic, for instance in \cite{DISTOLFO2016101,FRIES2011,HUO2020,HUANG2022,Aulisa2019} where new approaches to handle hanging nodes are proposed. In \cite{FRIES2011}, the methodology is explored in the context of the Extended Finite Element Method (G/XFEM), and in \cite{HUANG2022} the authors propose a new approach to handle hanging nodes in the context of the Scaled Boundary Finite Element Method (SBFEM). It is important to note that methodologies to refine the mesh while avoiding the creation of hanging nodes have been proposed in the literature, as for instance in \cite{Zienkiewicz2022,Rivara1984_1,Rivara1984_2,Rivara1992,Zhao2022}. The refinement patterns proposed in this chapter are an efficient alternative that can be used to refine meshes in optimal ways. The concept of refinement patterns has been used with different objectives as in \cite{Jaillet2022}, where they are associated with transition meshes between coarse and fine regions of a finite element mesh.

In this contribution, we focus on extending h adaptivity with refinement patterns. Refinement patterns are a powerful tool that provides a systematic way to handle complex refinements in regions of interest. The main idea is to generate a library of possible patterns on which elements can be refined and use this library to inform an h adaptive code on how to tackle complex cases. It is noted that even if a refinement pattern does not exist apriori, it can be readily generated given the nodes and sub-elements that compose a certain refinement pattern. Refinement patterns can be used to refine meshes in optimal directions as explained in Section \ref{sec:ref-pat-tools}. By integrating refinement patterns into an h adaptivity framework, one can improve the overall solution quality and computational efficiency. The methodology here explained has been integrated into the \texttt{NeoPZ} object-oriented environment for Finite Element (FE) analyses \cite{Devloo1997} and its implementation is used throughout the text to exemplify a possible way to integrate refinement patterns into an existing FE code.

The main contributions in this chapter are the following:
\begin{itemize}
  \item The concept of refinement patterns and their role in addressing hanging nodes in adaptive h-refinement are introduced.
  \item A detailed methodology for incorporating refinement patterns into the Finite Element Method is presented.
  \item The generality of the proposed approach is illustrated through mesh refinement of problems with complex geometries.
\end{itemize}

In this contribution, we intentionally do not include simulation results related to the application of refinement patterns. Our main objective is to present the methodology and concept behind it, rather than focusing on specific simulations.

This chapter is organized as follows. Sections \ref{sec:mat-descript} and \ref{sec:methodology} provide an overview of the fundamental concepts of refinement patterns and implementation aspects. Section \ref{sec:ref-pat-tools} presents necessary operations for the use of refinement patterns. Section \ref{sec:ref-pat-db} describes the methodology to generate a refinement pattern entry in the \texttt{NeoPZ} library. Section \ref{sec:Examples} presents complex refinement cases and discussions. Finally, Section \ref{sec:Conclusion} concludes the chapter and the outlines for future research directions are presented in Section \ref{sec:FutureWork}.

\section{Mathematical description \label{sec:mat-descript}}

The description of refinement patterns is based on the mathematical
concept of parametrization of the master elements and their sides.
These parametric spaces allow the formal description of the projection of
points from the interior of an element to the sides of the element
and vice-versa. They also allow to transfer point on an element/side
to a corresponding point on its ancestor/side. Finally, the transformation
between parametric spaces allows to formally describe shape function
restraints for $\hone$, $\hdivpure$ or $\hcurl$ spaces.

\subsection{Affine tranformations}

An affine transformation $T:R^{n}\rightarrow R^{m}$ has the form
\begin{align*}
x_{m} & =T(x_{n})\\
 & =A^{mn}x_{n}+B^{m},
\end{align*}
where $A^{mn}$ is an $m\times n$ matrix and $B^{m}$ an $m$-dimensional
vector. Affine transformations are utilized to describe the relationship between parametric spaces of different entities in a finite element mesh.

A projection $P:R^{n}\rightarrow R^{n}$ is an affine transformation
such that 
\[
P(x_{n})=P(P(x_{n})).
\]

When $m=3$ and $n=3$, affine transformations describe graphical transformations
commonly used in computer graphics. These transformations are so important
that they motivated the creation of the graphical processor unit (GPU) \cite{wikiGPU}.

\subsection{Element topology}

The element topology defines the parametric space associated with a
master element. Table \ref{tab:Element-topologies} enumerates the
element topologies included in the \texttt{NeoPZ} environment. As mentioned in Section \ref{sec:intro}, the \texttt{NeoPZ} environment is implemented within an object-oriented paradigm and, therefore, concepts can be "encapsulated" with the use of classes. In this context, a \texttt{class} is dedicated to each topology available in the environment.

\begin{table}[!ht]
  \centering
  \small
  \begin{tabular}{>{\centering\arraybackslash} m{2.1cm} >{\centering\arraybackslash} m{4.3cm} >{\centering\arraybackslash} m{3.7cm} >{\centering\arraybackslash} m{1.7cm} >{\centering\arraybackslash} m{3cm}}
    \hline
    \textbf{Name} & \textbf{Graphical representation} & \textbf{Parametric coordinates}                         & \textbf{Number of sides} & \textbf{Class name} \tabularnewline
    \hline
    Point         & \includegraphics*{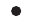}   & $\emptyset$                                                 & 1                        & \texttt{TPZPoint}  \tabularnewline
    \hline
    Line          & \includegraphics*{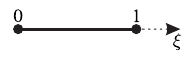}    & $-1\leq\xi\leq1$                                        & 3                        & \texttt{TPZLine} \tabularnewline
    \hline
    Quadrilateral & \includegraphics*{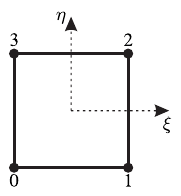}    & $-1\leq\xi\leq1$, $-1\leq\eta\leq1$                     & 9                        & \texttt{TPZQuadrilateral}  \tabularnewline
    \hline
    Triangle      & \includegraphics*{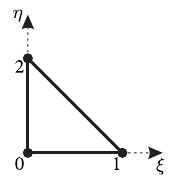}     & $0\leq\xi\leq1$, $0\leq\eta\leq1-\xi$                       & 7                        & \texttt{TPZTriangle}  \tabularnewline
    \hline
    Hexahedra     & \includegraphics*{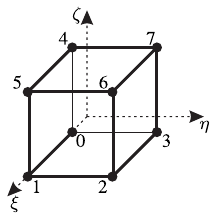}    & $-1\leq\xi\leq1$, $-1\leq\eta\leq1$, $-1\leq\zeta\leq1$ & 27                       & \texttt{TPZCube}  \tabularnewline
    \hline
    Tetrahedra    & \includegraphics*{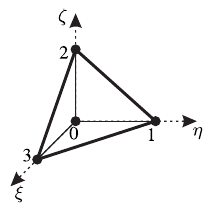}     & $0\leq\xi\leq1$, $0\leq\eta\leq1-\xi$, $0\leq\zeta\leq1-\xi-\eta$    & 15                       & \texttt{TPZTetrahedron} \tabularnewline
    \hline
    Prism         & \includegraphics*{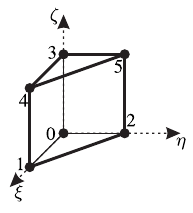}   & $0\leq\xi\leq1$, $0\leq\eta\leq1-\xi$, $-1\leq\zeta\leq1$   & 21                       & \texttt{TPZPrism} \tabularnewline
    \hline
    Pyramid       & \includegraphics*{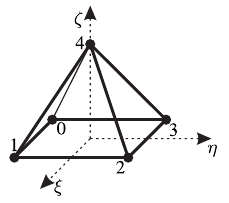}     & $-1+\zeta\leq\xi\leq1-\zeta$, $-1+\zeta\leq\eta\leq1-\zeta$, $0\leq\zeta\leq1$  & 19                       & \texttt{TPZPyramid}  \tabularnewline
    \hline
  \end{tabular}
  \caption{Element topologies in local (or master) coordinate system.}
  \label{tab:Element-topologies}
\end{table}

\begin{figure}[!ht]
  \centering
  \includegraphics{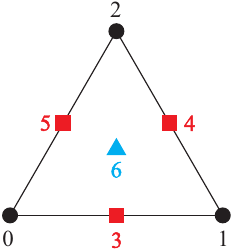}
\caption{Topology of a triangular element composed of seven sides: three zero-dimensional sides (points) with local indices 0, 1, and 2; three one-dimensional sides (edges) with local indices 3, 4 and 5; and one two-dimensional side (face) with local index 6. \label{fig:Triangle-topology}}

\end{figure}

Each master element can be viewed as the union of open sets of points
associated with sides. Each side is in turn associated with a topology
and has an associated parametric space. The definition of the sides
and their orientation is an integral part of the definition of the
topology. Figure \ref{fig:Triangle-topology} illustrates the master
element of a triangle. Note that the triangle is composed of 3 points (zero-dimensional sides),
3 edges (one-dimensional sides) and one face (two-dimensional side). An affine transformation
can be defined between the parametric spaces associated with the sides
and the master element space.

\subsubsection*{Definitions}

The following definitions are used in the description of the topology
\begin{itemize}
\item $T_{es}$: Affine transformation that projects a point from the
interior of the element to the parametric space of side $s$.
\item $T_{se}$: Affine transformation transforming the location of
a point defined in a parametric space of side $s$ to a point in the
element space.
\end{itemize}
 
\subsubsection*{Consistency checks}

If a point $x$ is such that it lies on a side then, there exists a
point in parametric space of $s$, denoted $x_{s}$, such that
\[
x=T_{se}x_{s}.
\]

Also, if $T_{es}$ projects a point in the element space onto a side
and assuming that $x$ lies on the side, then
\[
x_{s}=T_{es}x,
\]
implying that
\[
x_{s}=T_{es}T_{se}x_{s}\;\forall x_{s}.
\]

Therefore, taking the composition of transforming a point in the parametric
space of a side $s$ to the element space and transforming it back
to the side results in an identity operator (which is in fact a
projection)

\[
I_{s}=T_{es}T_{se}.
\]

Defining the affine transformation $P_{s}$ as the composition:

\[
P_{s}=T_{es}T_{se},
\]
then $P_{s}$ is a projection operator such that 
\[
P_{s}P_{s}\equiv P_{s}.
\]

Then clearly
\begin{align*}
P_{s}P_{s} & =T_{es}T_{se}T_{es}T_{se}\\
 & =T_{es}I_{s}T_{se}\\
 & =P_{s},
\end{align*}
meaning that if a point from the element space is projected to the
parametric space of a side and then back to element space, then this is in fact a projection.

\subsubsection*{Composition of projections}

If the closure of a  face $f_{a}$ of a three-dimensional element contains the one-dimensional side $s_{c}$ then the composition of the projection to the face $P_{f_a}$ and and the projection to side $s_{c}$ (i.e. $P_{s_c}$) is equal to the projection to $s_{c}$
\[
P_{s_{c}}=P_{s_{c}}P_{f_{a}}.
\]

The description of the topology is the common denominator that relates the description of the geometry, the definition of the approximation space, and the definition of a numerical integration rule. 

Each master element defines its own set of sides and orientations.
In the \texttt{NeoPZ} library the element topologies are defined in the \texttt{Topology}
subdirectory of the library. Each topology \texttt{class} is \texttt{static}, meaning
its behavior does not depend on an internal data structure. The parametric
transformations between sides for the different topologies are documented in \cite{DEVLOO2009}.

\subsection{Neighboring information \label{sec:neigh-info}}

\begin{figure}[!ht]
  \centering
  \includegraphics{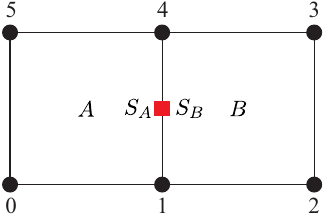}
  \caption{Illustration of two neighboring elements through the side with global node indices 1 and 4.\label{fig:Illustration-of-element}}
\end{figure}

The concept of viewing the master element as the union of open sets
of points allows to define the concept of neighbors. Let $A$ and $B$ be two elements of a finite element mesh and $S_A$ and $S_B$ sides of these two elements. Then, element $A$ is the neighbor of element $B$ through side $S_A$ if and only if the set of points associated with $A/S_A$ is equal to the set of points
associated with $B/S_B$. Figure \ref{fig:Illustration-of-element}
illustrates the concept of neighbors graphically.

\begin{figure}[!ht]
  \centering
  \includegraphics{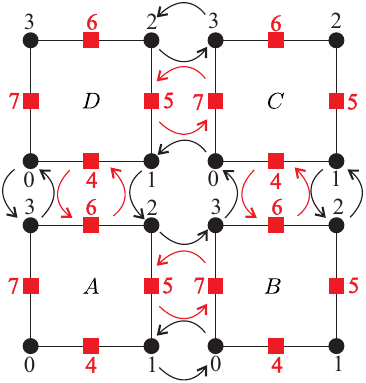}
\caption{Circular data structure of neighboring elements. Black represents neighboring data structures between nodes and red represents neighboring data structures between edges. \label{fig:Element-neighbors}}

\end{figure}

An element $A$ along the side $S_A$ can be a neighbor to more than one element/side.
Dynamic memory allocation is avoided by keeping track of element neighbors
in a circular structure. Figure \ref{fig:Element-neighbors} illustrates
a point that is shared by 4 elements. The circular data structure represents
a linked list:
\[
A/S_A\rightarrow B/S_B\rightarrow C/S_C\rightarrow D/S_D\rightarrow A/S_A
\]

As such the data structure associated with side $S_A$ of element $A$
only keeps track of element $B/S_B$. By cycling through the neighboring
information of the element/sides, one can collect the complete set
of elements/neighbors. As the number of sides of $A$ is fixed (associated
with the topology), a fixed-sized data structure keeps track of the
neighboring information.

\subsubsection*{Building the neighboring information}

Considering that the data structure of neighboring information is
not standard, we describe a simple algorithm to build the neighboring
information based on a \emph{traditional} element node data structure.

The first premise is that each side of an element can be associated
with a set of nodes. The closure of the side contains a
given number of nodes. The global index of these nodes forms the set
of nodes associated with a side.

The second premise is that, in a traditional finite element mesh, two
elements are neighbors along one side if they share the associated set
of nodes, \emph{i.e.} the side is uniquely determined by the nodes (\emph{e.g.}
if two triangular sides share three nodes, then they are neighbors
along a face).

\subsubsection*{$\mathbf{A}$ Neighbors along nodes}

The first step is to loop over the elements once to fill in the neighboring
information along the nodes. To do so, a temporary vector consisting of pairs of integers of the size of the number of nodes in the mesh is created. This data structure is denoted \texttt{NeighInfo} and will contain at each position one possible (element index)/(side index) pair associated with the global index of a node. Then the algorithm to fill in the neighboring information along nodes is shown in Algorithm \ref{alg:neigh-nodes}.

\begin{algorithm}[!h]
  \SetKwInOut{Input}{Input}
  \SetKw{Continue}{continue}
  \SetKw{Break}{break}
  \Input{Finite element mesh and \texttt{NeighInfo}}
  \For {Element $e^i$ in the Finite Element Mesh}
  {
    \For {Node $n^j$ of $e^i$}
    {
      \If{Position $k$ of $\texttt{NeighInfo}$ is not initialized where $k$ is the global index of $n^j$}
      {
        Insert the pair $e^i$/$n^j$ in position $k$ of \texttt{NeighInfo}
      }
      \Else
      {
        $e_i$/$n_j$ is inserted in the neighboring data structure of the element contained in position $k$ of \texttt{NeighInfo}
      }
    }
  }
  \caption{Algorithm to compute all element neighbors through nodes. \label{alg:neigh-nodes}}
\end{algorithm}

\subsubsection*{$\mathbf{B}$ Neighbors along higher dimensional sides}

Let \texttt{NeighSideInfo} be a temporary data structure created for each side of an element of size equal to the number of nodes of that side. Each position in \texttt{NeighSideInfo} pertains to a certain node of the side and will contain a set of integers of all the neighboring elements through that node. Note that the node neighboring data structure for each element has been previously computed in Algorithm \ref{alg:neigh-nodes}. Algorithm \ref{alg:neigh-sides} is used to compute all element neighbors through sides of dimension higher than 0.

\begin{algorithm}[!h]
  \SetKwInOut{Input}{Input}
  \SetKw{Continue}{continue}
  \SetKw{Break}{break}
  \Input{Finite Element mesh}
  \For {Element $e^i$ in the Finite Element Mesh}
  {
    \For {Side $s^j$ of $e^i$}
    {
      \If{dimension of $s^j$ is 0}
      {
        \Continue
      }
      Create \texttt{NeighSideInfo} of size equal to the number of nodes of $s^j$

      Fill \texttt{NeighSideInfo} based on data structure of neighboring nodes created using Algorithm \ref{alg:neigh-nodes}

      Perform intersection of the sets in \texttt{NeighSideInfo}

      The result is the list of elements that are neighbors of $e^i$ through side $s^j$
    }
  }
  \caption{Algorithm to compute all element neighbors through sides of dimension higher than 0. \label{alg:neigh-sides}}
\end{algorithm}

It is noted that the complexity of this algorithm is linear. It can possibly be parallelized.

\subsection{Parametric transform between neighbors}

Two elements are neighbors if they share a side. The parametric transform between both elements depends on the relative orientation of both sides. The number of possible combinations depends on the topology of the side. A line has two possible orientations, a triangle has six possible orientations, and a quadrilateral has eight possible orientations. 

The correctness of the parametric transformation can verified by the following tests: for an arbitrary point on side $S_A$ of element $A$. The corresponding point on side $S_B$ of element $B$ is computed by the transformation $T_{AB}$. The $x$ coordinate computed by element $A$ should correspond to the $x$ coordinate computed by element $B$. 

\subsection{Parametric transform between son side and father side}

It is an essential feature of geometric refinements that a {\em father} element is partitioned by its son elements. Moreover, each side of the son elements is entirely contained in a side of the father element. As a consequence of this property, an affine transformation can be defined between the parametric space of the son/side and the corresponding father/side. This affine transformation is computed by an $L^2$ projection. Let $A^{mn}$ and $B^{m}$ be the matrix and vector of the affine transformation, respectively. The matrix $A^{mn}$ and vector $B^m$ are computed by minimizing the following functional:

\begin{equation}
\min_{A^{mn},B^m} \int_{\hat{S}} ||x_S(\xi) - x_F(A^{mn}\xi - B^m)||^2 d\xi,
\end{equation}
where $\hat{S}$ is the parametric space of the son/side, $x_S$ is the cartesian coordinate of the son/side, and $x_F$ is the cartesian coordinate of the father/side. In Section \ref{sec:sub-el-param-transf}, the same functional is used to compute the affine transformation between the parametric space of the father element and the parametric space of the son element within the context of the \texttt{NeoPZ} environment.

\section{Data structure and implementation aspects \label{sec:methodology}}

In the following subsections, details regarding a choice of implementation of refinement patterns in the \texttt{NeoPZ} library \cite{Devloo1997} are discussed. In Subsection \ref{sec:fund-ref-pat}, preliminary concepts are introduced. Subsection \ref{sec:def-ref-pat} describes
the data structure required for defining a valid refinement pattern and Subsection \ref{sec:ini-ref-pat}
discusses the operations to load refinement patterns into an FE code. Finally, in subsection \ref{sec:ref-pat-tools}, useful routines are presented in a \emph{bottom-up} fashion,
illustrating the building blocks needed to define the routine \texttt{RefineDirectional} and providing
insights on the types of operations possible using these patterns. For more details on the implementation, the reader is referred to the classes \texttt{TPZRefPattern} and \texttt{TPZRefPatternTools} in the \texttt{NeoPZ} library available on \href{https://github.com/labmec/neopz}{Github}.

\subsection{Fundamental concepts of refinement patterns}\label{sec:fund-ref-pat}
The term refinement pattern is used here for describing a systematic way of geometrically refining an element, \emph{i.e.}, splitting it into more elements. This procedure is performed in the reference (or master) element and, therefore, a given refinement pattern can be readily applied to different elements of the same type. In a given mesh, each element may also have refinement patterns associated with it because of \emph{compatibility} with the refinement patterns of its neighboring elements. This means that the vertices created by the refinement pattern of one element must be compatible with the vertices created by the refinement pattern of its neighbors.

In practice, to reduce the number of needed refinement patterns, one idea is to use the concept of \emph{permutation} of refinement patterns: since the refinement pattern is defined in the reference element, changing the node ordering of such element will result in different patterns. Being able to compute the permutations of a refinement pattern allows one to apply them to elements with a different node ordering.

It is noted that dividing an element into subelements using refinement patterns can either improve or worsen the aspect ratio of the original element. Using refinement patterns neither ensures nor impedes the improvement of the mesh quality.

\subsection{File structure defining a refinement pattern}\label{sec:def-ref-pat}

\begin{figure}
\begin{minipage}{0.5\textwidth}
  \begin{codingtime}
    \small
    \textcolor{gray}{\% int(\# nodes) int(\# elements)}
    
    7 3
    
    \textcolor{gray}{\% int(id) string(name)}
    
    -50 TetraDir1Side 

    \bigskip
    \textcolor{gray}{\% double[3](node coords)}
    
    0 0 0 

    1 0 0

    0 1 0

    0 0 1

    0 0 0.5

    0.5 0 0.5

    0 0.5 0.5

    \bigskip
    \textcolor{gray}{\% int(el type) int(mat id) int(node ids)}
    
    4 1 0  1  2  3 
    
    6 1 0  1  2  4  5  6
    
    4 1 4  5  6  3
  \end{codingtime}
\end{minipage}
\begin{minipage}{0.5\textwidth}
  \centering
  \includegraphics[width=\textwidth]{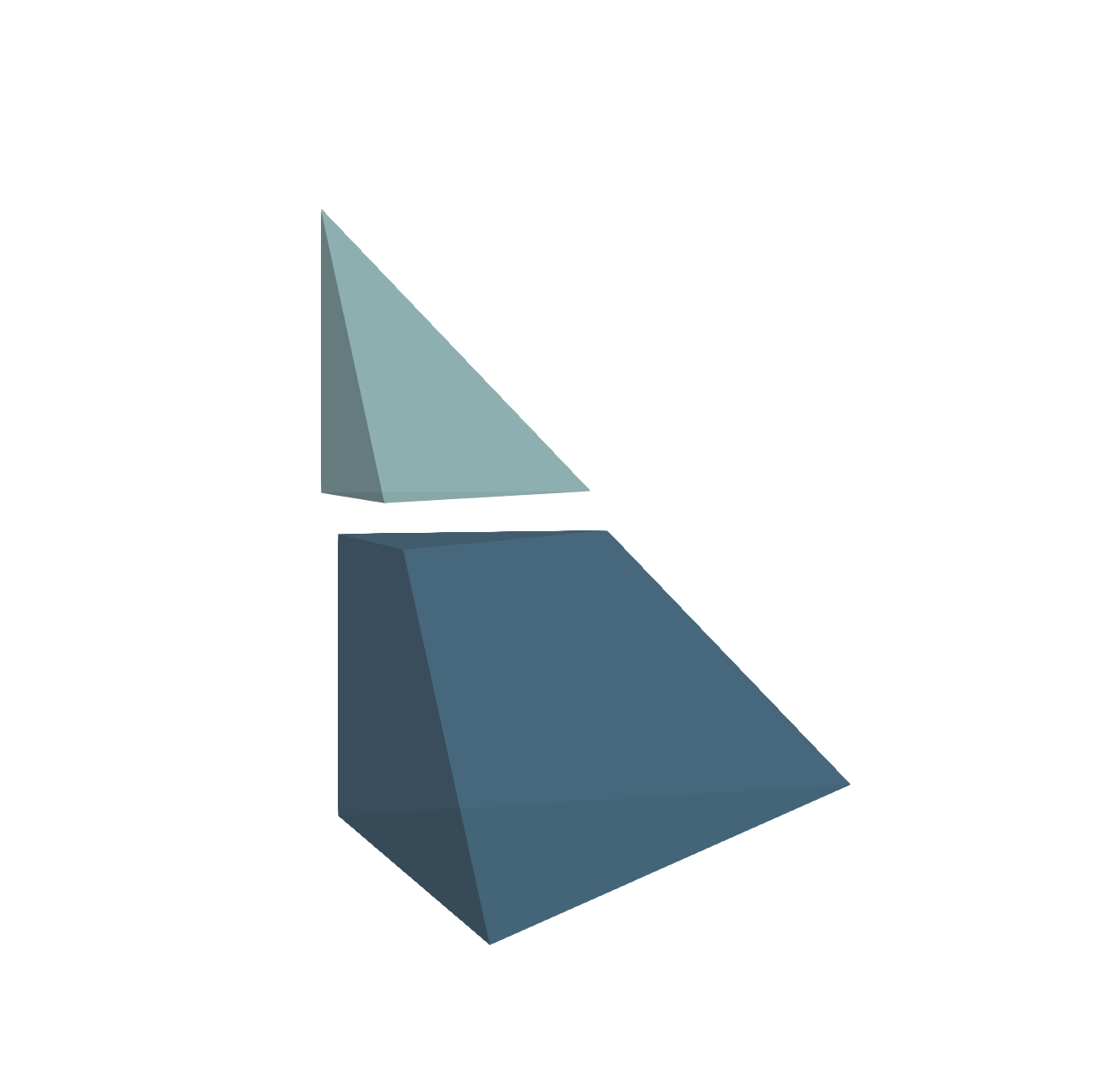}
\end{minipage}
\caption{Description of a simple refinement pattern on a tetrahedron and the resulting elements.}\label{fig:ref-pat-algo}
\end{figure}

A refinement pattern is described through a secondary geometric mesh -- different than the finite element geometric mesh used for the simulation -- whose domain $\Omega$ corresponds to a given
supported element type (linear, triangular, quadrilateral, hexahedral, \emph{etc.}). It can thus be described in text format as in Figure \ref{fig:ref-pat-algo}.

Figure \ref{fig:ref-pat-algo} illustrates a tetrahedron split into another tetrahedron and a wedge. The data file is divided into three sections. In the first section, we specify the number of nodes (\texttt{7}), the number of elements (\texttt{3}, including the original element that is refined), and provide both a numerical identifier and a name that are associated with the pattern. The second section consists of the numerical coordinates of the nodes in the three-dimensional domain. The third section specifies all the elements in the mesh, possibly assigning different identifiers for them. The first element to be described is always the father element and the subsequent are the children. In this example, the tetrahedral (hexahedral) element is represented by the integer \texttt{4} (\texttt{6}), and all the elements are assigned the same identifier (\texttt{mat id}), \texttt{1}.

Once a given refinement pattern is read using the file described in this section, the next step is to load it into the memory and use it to perform the necessary operations for applying this pattern to refine an element.

\subsection{Creating a refinement pattern}\label{sec:ini-ref-pat}

As discussed in the previous section, a refinement pattern consists of a secondary geometric mesh corresponding to a partition of a master element $K$. In this section, the procedures performed on this secondary mesh to apply a refinement pattern to a given element are discussed.

\subsubsection{Reading the pattern}
Once this secondary geometric mesh is read, the first operation is to convert the coordinates of the nodes in $K$ to the ones in $\hat{K}$, the \emph{reference element}. If the original coordinates do not correspond to the coordinates of the master element, the inverse of the map $F: \hat{K}\rightarrow K$ is used. This is important since refinement patterns are invariant under affine transformations.

Next, the usual connectivities of the mesh are computed, as described in Section \ref{sec:neigh-info}. Then, the transformation between the possible \emph{permutations} of the element is computed, as described in Subsection \ref{sec:fund-ref-pat}. By doing so, the pattern is available regardless of the node ordering of the element in the mesh. After these procedures, refinement-specific tasks have to be performed. 

\subsubsection{Sub-element parametric transformations}\label{sec:sub-el-param-transf}
To obtain conforming approximation spaces, each sub-element in the refinement pattern must be aware of
its relationship with the parent element so that the relevant traces of basis functions can be compatibilized. More specifically, each side of each sub-element needs to know which side in the parent element contains it and what is the parametric transformation to that side in the parent element in parametric coordinates.

\begin{algorithm}
  \SetKwInOut{Input}{Input}
  \Input{sub-elements $K^e$ and parent element $K^f$}

  \For{ sub-element $K^e\in\{K\}$ }
  {
    \For {side $s^e\in \partial K^e$}
    {
      Compute center of side $s^e$ in $\xi_f$ in the parametric coordinates of $K^f$
      
      Find which side $s^f$ of the father contains the center point of $s^e$

      Compute the affine transformation $F_e^f: \xi_{s^e} \rightarrow \xi_{s^f}$
      
      Store $(s^f, F_e^f)$
    }
  }
  \caption{Computing parametric transform from the sides of the sub-elements to the side of the father element.}
\end{algorithm}

The transformation $F_e^f$, mapping a point in the parametric coordinates of the sub-element side
$s^e$ to a point in the parametric coordinates of the parent element's side $s^f$ is computed by
a minimization procedure. Namely, we search for the affine transform $F_e^f: A\xi_{s^e}+B$ by
finding the coefficients of $A$ and $B$ that minimize

\begin{equation}      
  \lVert\Delta\rVert^2=\int_{\Omega_{s^e}}(\xi_{s^f}-(A\xi_{s^e}+B))\cdot(\xi_{s^f}-(A\xi_{s^e}+B))\,d\xi_{s^e}\text{.}
\end{equation}

It is worth noting that the resulting $F_e^f$ is exact only if the jacobian of the mapping
of the sub-element is constant.

\subsubsection{Parent element side partitioning}

For each side of the father element, we store its \emph{partition}. By partition, we mean a data structure containing:

\begin{itemize}
\item Indices of internal nodes that result from the refinement pattern that is contained in it.
\item List of pairs of sub-element/side contained in it.
\end{itemize}

The partition of sides is used to compatibilize refinement patterns. This procedure can be better explained with an example. Assume that an element in a given mesh needs to be refined with a certain refinement pattern. Depending on the refinement pattern, new nodes are created in the mesh, that are contained in sides of the original element. The proposed paradigm dictates that these nodes should match nodes created on the corresponding side of a previously refined neighboring element, \emph{i.e.}, their side-refinement patterns should be compatible (Side refinement patterns are the subject of the next section). However, this verification can be costly. Storing the partition can be used to check that a given pattern is not compatible significantly faster.

This algorithm can be easily computed with the data from Subsection \ref{sec:sub-el-param-transf}, so it is omitted.

\subsubsection{Side refinement patterns}

A refinement pattern of an element implicitly defines refinement patterns for its sides. The refinement pattern of a side is defined by the refinement pattern of the element and the parametric transformation between the side and the element. When a neighboring element is refined, the refinement patterns associated with the neighboring sides should be identical. In other words, two neighboring elements have compatible refinement patterns if they have the same refinement pattern associated with their common side.

Embedded in a refinement pattern are the refinement patterns associated with each side of the element\footnote{Except, obviously, for zero-dimensional sides.}, and they are used in operations associated with boundaries and interfaces between elements.

The creation of the side refinement patterns is performed recursively: the mesh associated with the side is created and its pattern is initialized (if not found in a database of previously loaded refinement patterns, as discussed in Section \ref{sec:ref-pat-db}).

When a refinement pattern is applied to an element, the {\em code} checks whether the neighboring elements have compatible side refinement patterns. If they do not, the code will flag an error and abort. To assist the user in using the library, the code will also provide a list of compatible refinement patterns.

\section{Refinement pattern tools}\label{sec:ref-pat-tools}

With the extension of the \texttt{NeoPZ} library to include refinement patterns, a set of tools is necessary to perform operations on these patterns: whereas refinements using uniform refinement patterns are always compatible, this is not the case when using general refinement patterns. The tools presented in this section are used to verify the compatibility of refinement patterns and to find compatible refinement patterns that match certain criteria.

\subsection{RefPatternEquality}
Two refinement patterns are said to be equal if, and only if:
\begin{itemize}
\item For each element type, the same number of elements is present in both patterns
\item The number of nodes in both patterns is the same
\item There is a bijection between the nodes on both patterns
\item Every sub-element has the same (mapped) nodes
\end{itemize}

\subsection{GetCompatibleRefPatterns}
The refinement patterns are said to be compatible if, for every neighbor with a refined side, the element's side is refined exactly the same way. That implies that if no neighbors are refined, every refinement pattern is said to be compatible.
Thus, compatibility refers to the fact that applying a given pattern on a certain element in a mesh will result in a conforming mesh. This procedure is described in Algorithm \ref{alg:ref-pat-tools}.

\begin{algorithm}
  \SetKwInOut{Input}{Input}
  \SetKw{Continue}{continue}
  \SetKw{Break}{break}
  \Input{element $K_e$}
  \For {side $s^i\in \partial K_e$ such that $\dim{s^i}>0$}
  {
    \For {neighbor $K_n$ of ($K,s^i$)}
    {
      \If{There is side $s^n$ of $K_n$ that is refined}
      {
        Store refinement pattern\footnotemark of neighbors $P_n^{s^i}$
        
        \Break
      }
    }
  }
  \For {candidate pattern $P$ of $K_e$}
  {
    \For {side $s^i\in \partial K_e$ such that $\dim{s^i}>0$}
    {
      
      \If{$P_n^{s^i}$ exists}
      {
        \If{side pattern of $P$ is different from $P_n^{s^i}$}
        {
          Mark as non-compatible

          \Break
        }
      }
      
    }
    \If{$P$ is not marked as non-compatible}
    {
      Insert into the compatible list
    }
  }
  \caption{\texttt{GetCompatibleRefPatterns} \label{alg:ref-pat-tools}}
\end{algorithm}
\footnotetext{Appropriately permuted as to match node ordering of $s^i$}

\subsection{PerfectMatchRefPattern}
The \texttt{PerfectMatchRefPattern} method receives as input a list of marked edges and, among all the compatible refinement patterns returned from \texttt{GetCompatibleRefPatterns}, finds the patterns that do not refine any edge that is not marked. These are considered to be a \emph{perfectly matched}, meaning that a set of patterns satisfying certain criteria requiring the least amount of refinement has been found. Currently, in the \texttt{NeoPZ} library, the pattern with the least elements in the set of perfect match patterns is chosen.

\texttt{PerfectMatchRefPattern} corresponds to the selection of a refinement pattern that creates {\em transition} refinement patterns to create meshes without hanging nodes.

\subsection{RefineDirectional} 

Lastly, the \texttt{RefineDirectional} method is presented, which has extensive applicability in problems where graded meshes are desired in the direction of a surface, line or point. The idea of refining towards a certain direction has been explored in \cite{Rivara1984_2,Zienkiewicz1994} but limited to simplicial meshes. 

The main idea of this method is to refine elements that are neighbors to a given entity (point, line, surface) in such a way that the entity itself is not divided. All the edges that touch the entity and are not included in the entity are divided. The algorithm is described in Algorithm \ref{alg:ref-dir}.

\texttt{RefineDirectional} are illustrated in the examples of Section \ref{sec:Examples}.

\begin{algorithm}
  \SetKwInOut{Input}{Input}
  \SetKw{Continue}{continue}
  \SetKw{Break}{break}
  \Input{List of elements that can be refined $\{K\}$, identifier $m$ of the region of interest}
  \For{element $K\in \{K\}$}
  {
    \For{vertex $v$ of $K$}
    {
      \If{$(K,v)$ has a neighbor with identifier $m$}
      {
        Marked $v$
      }
    }
    \If{No vertex has been marked}
    {
      \Return
    }
    \For{edge $e$ of $K$}
    {
      \If{$(K,e)$ has a neighbor with identifier $m$}
      {
        \Continue
      }
      \If{only one vertex $v_e$ of $e$ is marked}
      {
        Mark edge $e$
      }
    }
    Find \text{PerfectMatchRefPattern} with marked edges
  }
  \caption{\texttt{RefineDirectional}}
  \label{alg:ref-dir}
\end{algorithm}

\section{Refinement pattern database}\label{sec:ref-pat-db}
All the operations performed so far depend on the availability of different refinement patterns at runtime. Therefore, one essential part of the usage of refinement patterns is to appropriately store them in a database, which, in the \texttt{NeoPZ} library is implemented as a globally accessed instance of the class \texttt{TPZRefPatternDataBase}.

The database is populated with refinement patterns that are read from files.
When reading a refinement pattern from a file, the \texttt{RefPatternEquality} algorithm of Section \ref{sec:ref-pat-tools} is used to check if this refinement pattern is already present in the database. If it is not, it is inserted, and indexed by the numerical identifier and name, as shown in Figure \ref{fig:ref-pat-algo}.

The indexing of the refinement patterns makes searching for a pattern in the database correspond to a map lookup operation, therefore with complexity $\mathcal{O}(\log n)$, where $n$ stands for the number of stored refinement patterns. 

For each refinement pattern read from a file, its permutations and side refinement patterns, if not in the database, are added to the database as well.

It is noted that the database also represents a crucial saving of memory for a finite element mesh supporting refinement patterns: each element in the mesh can store a pointer to its refinement pattern, instead of storing the pattern itself. Taking into account that permutations and side refinement patterns are also needed at the element level, it is clear that avoiding duplicate instances of patterns becomes essential when the number of available refinement patterns grows.

All the refinement patterns in the database are available to be used in routines such as \texttt{Get\-Compatible\-Ref\-Patterns}, which is called in the operations such as directional refinement (\texttt{Refine\-Directional}). Therefore, it is important to populate the database with useful patterns. This is done through a list of patterns stored in text files, which are read before making use of the routines of Section \ref{sec:ref-pat-tools}.

\section{Examples\label{sec:Examples}}

This section is devoted to demonstrating some of the capabilities of the refinement tools that were discussed in the previous sections.

\subsection{YF-17 aircraft mesh}

\begin{figure}[!ht]
  \centering
  \subfloat[]{\includegraphics[width=0.47\textwidth]{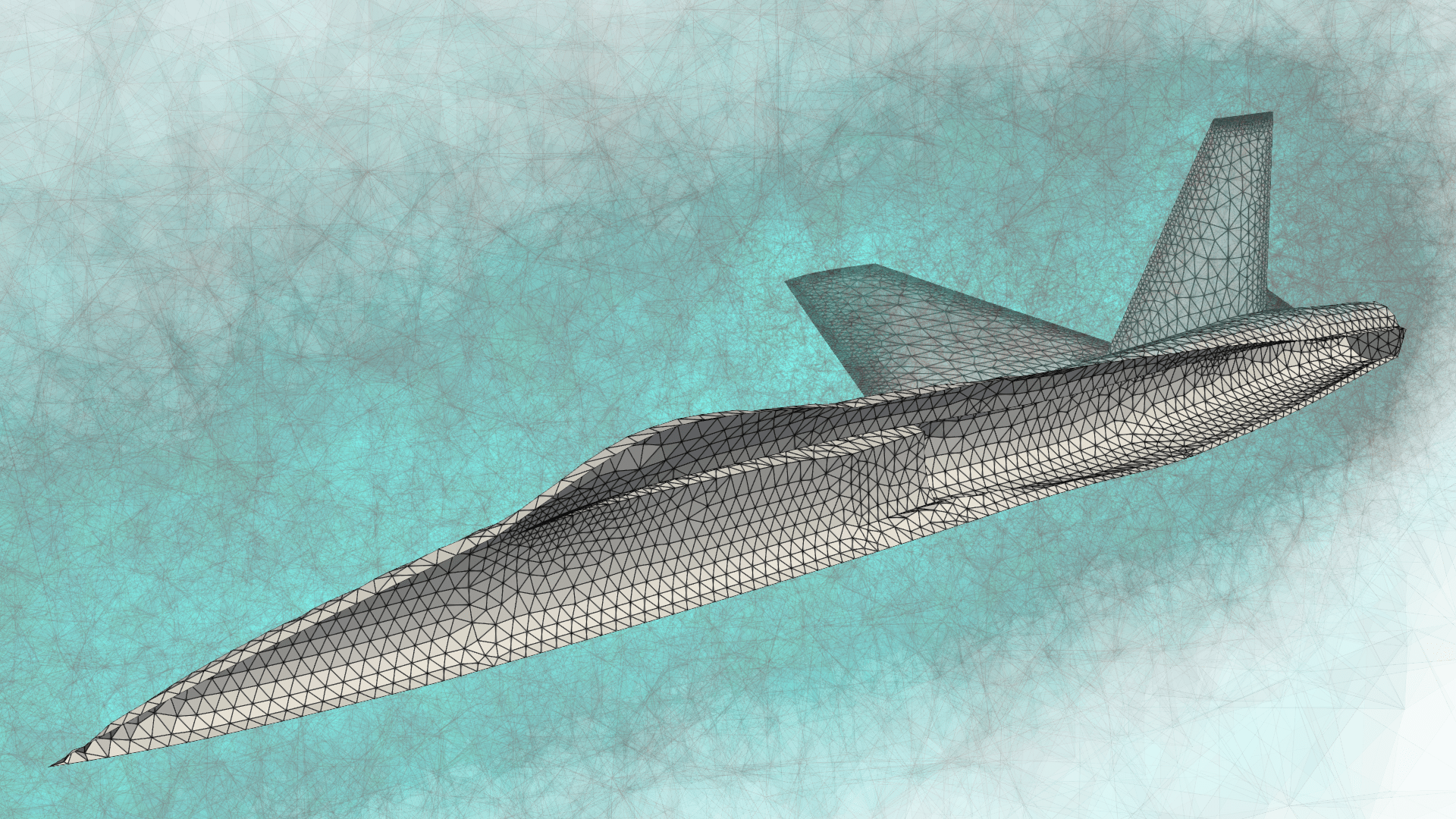}} \hfill
  \subfloat[]{\includegraphics[width=0.47\textwidth]{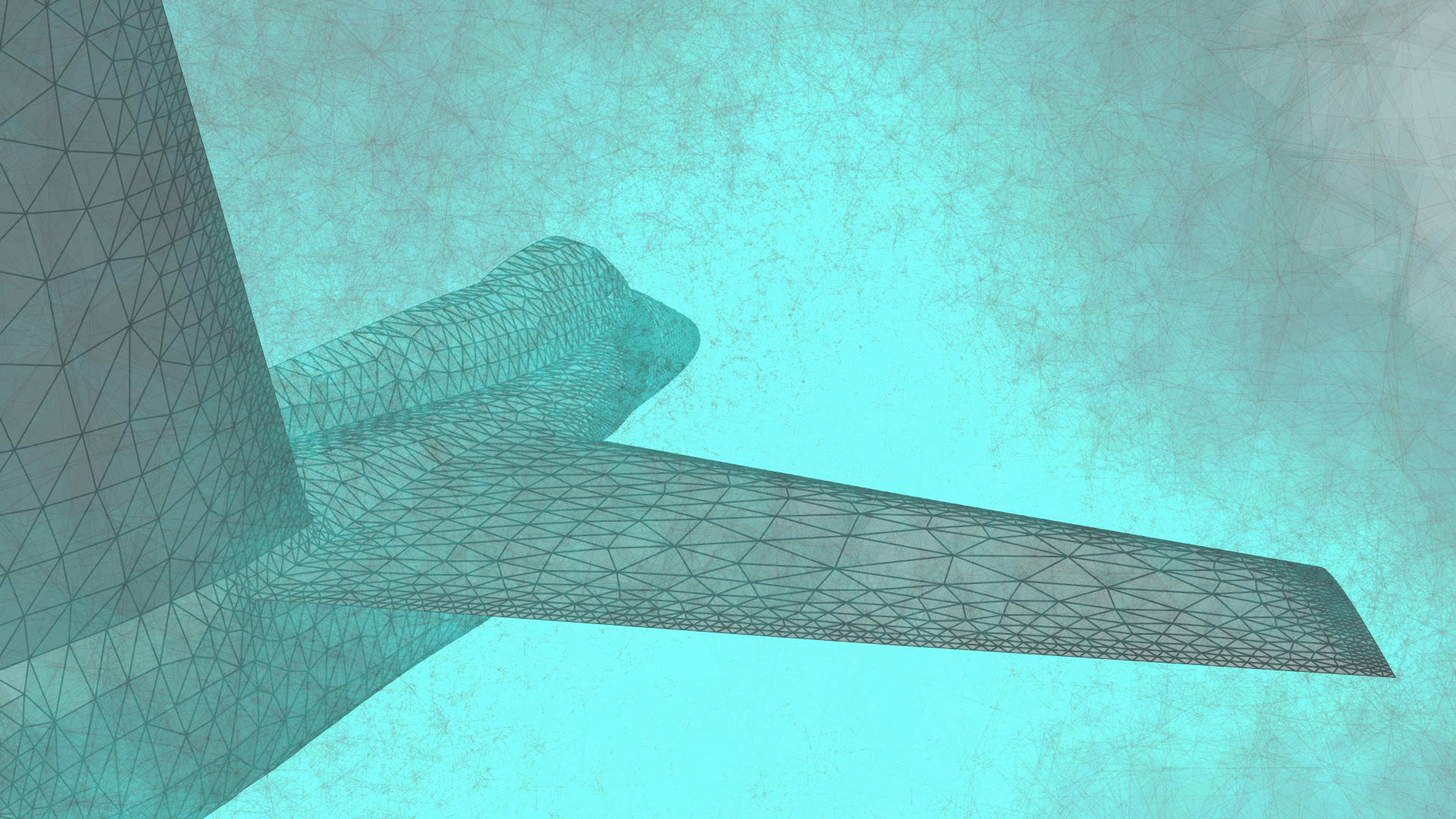}}
  \caption{YF-17 aircraft mesh. (a) Perspective view and (b) rear view.}
  \label{fig:YF17-surface}
\end{figure}

To illustrate the usage of \texttt{RefineDirectional} routine, we consider a surface of a YF-17 fighter jet, available in \cite{f17mesh}, embedded in a tetrahedral mesh that represents its surrounding air. The initial mesh is shown in Figure \ref{fig:YF17-surface}, and consists of 359\,696 tetrahedral elements. 

\begin{figure}[!ht]
  \centering
  \includegraphics{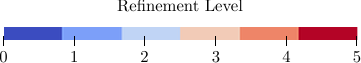}\\
  \subfloat[Level 0]{\includegraphics[width=0.49\textwidth, trim={0, 5cm, 0, 5cm}, clip]{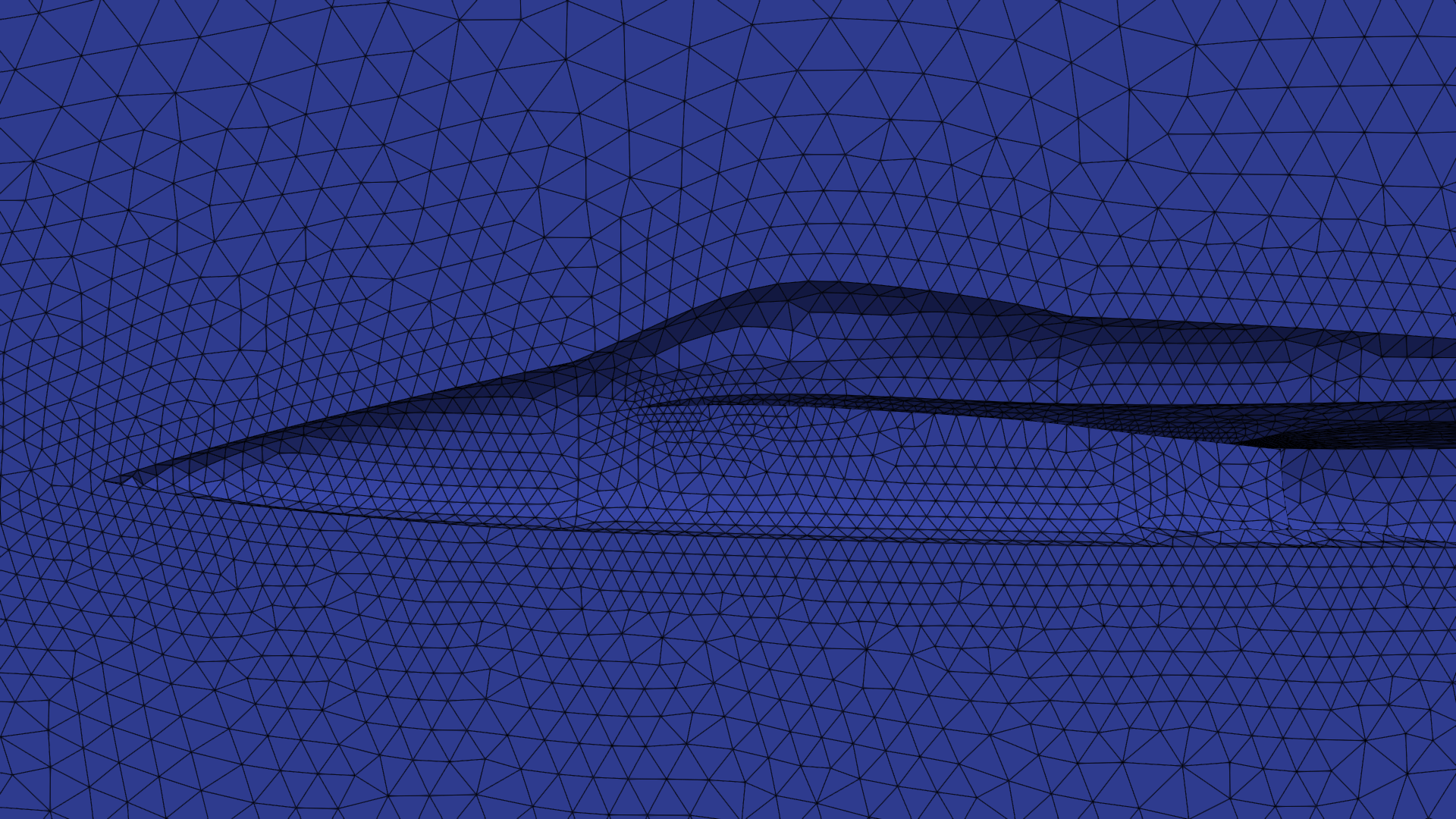}} \hfill
  \subfloat[Level 1]{\includegraphics[width=0.49\textwidth, trim={0, 5cm, 0, 5cm}, clip]{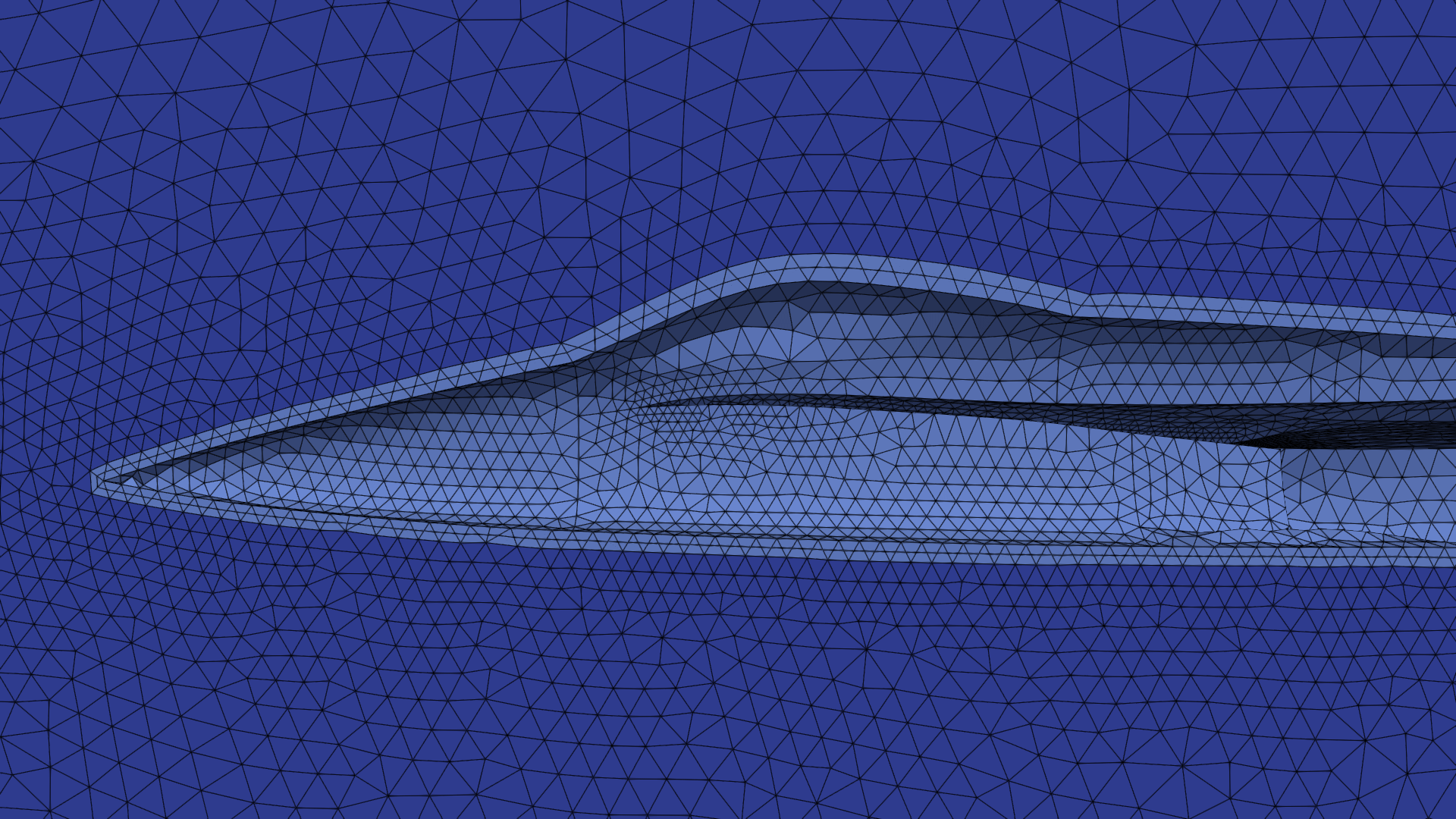}} \\
  \subfloat[Level 2]{\includegraphics[width=0.49\textwidth, trim={0, 5cm, 0, 5cm}, clip]{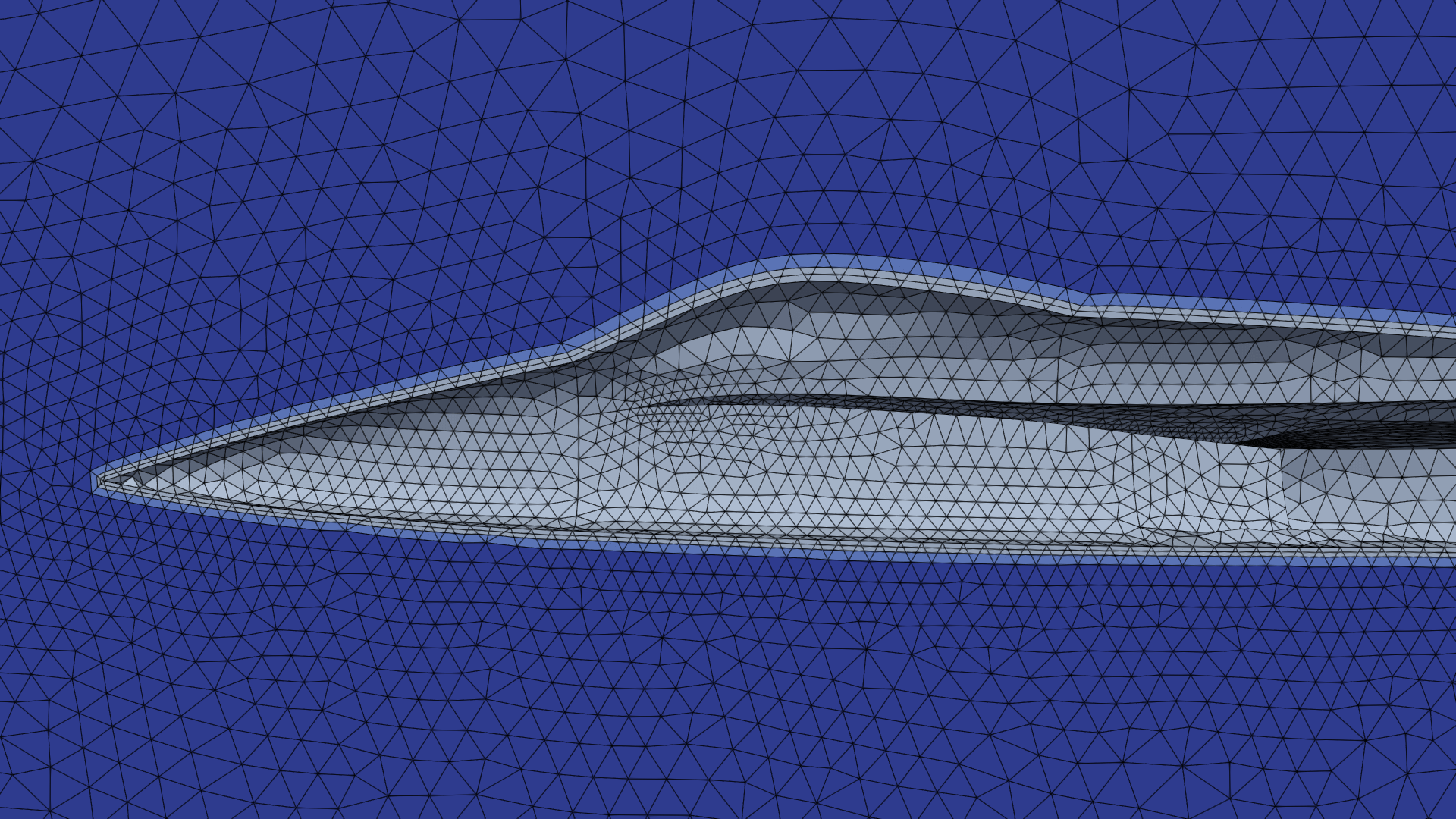}} \hfill
  \subfloat[Level 3]{\includegraphics[width=0.49\textwidth, trim={0, 5cm, 0, 5cm}, clip]{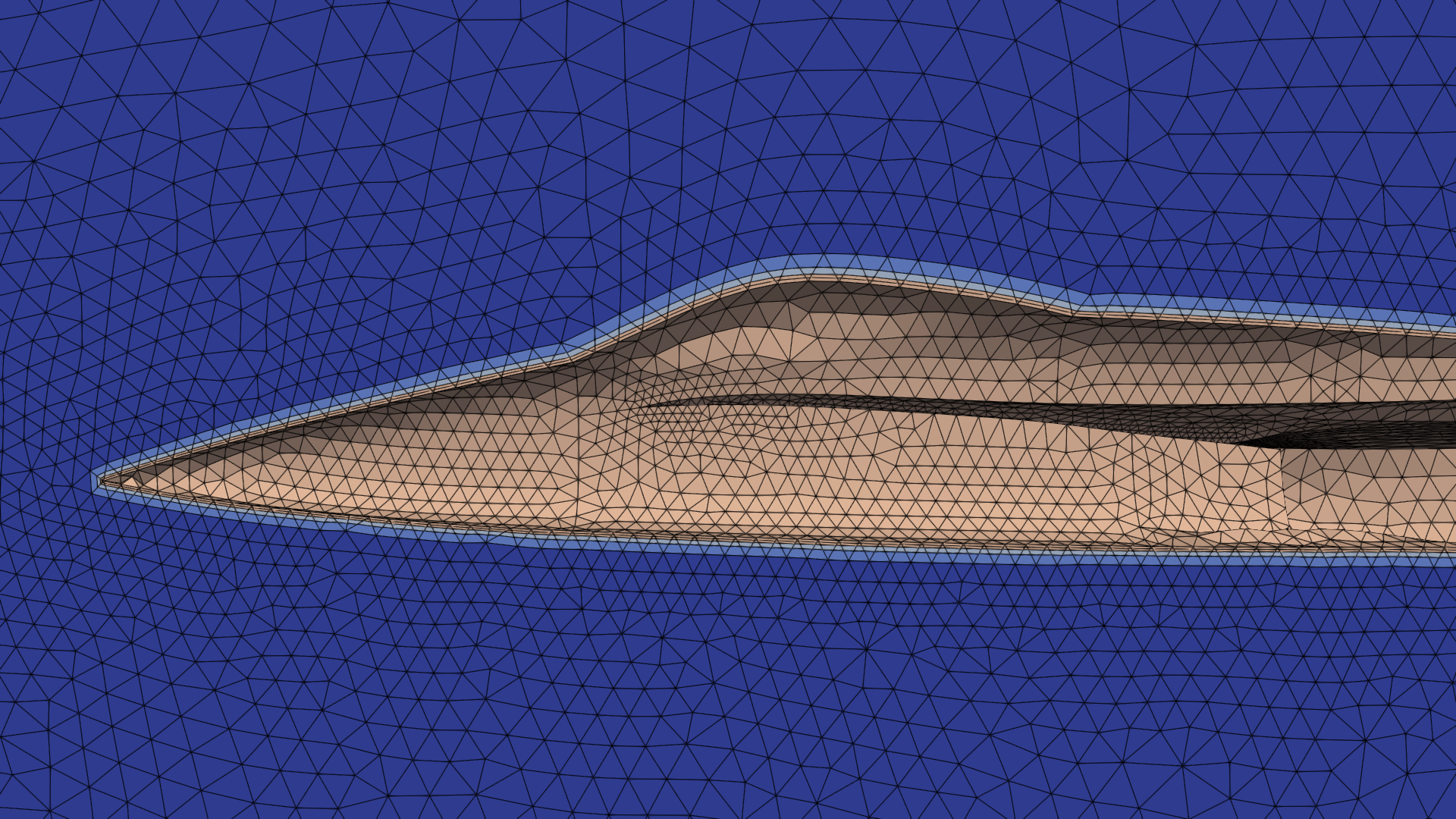}} \\
  \subfloat[Level 4]{\includegraphics[width=0.49\textwidth, trim={0, 5cm, 0, 5cm}, clip]{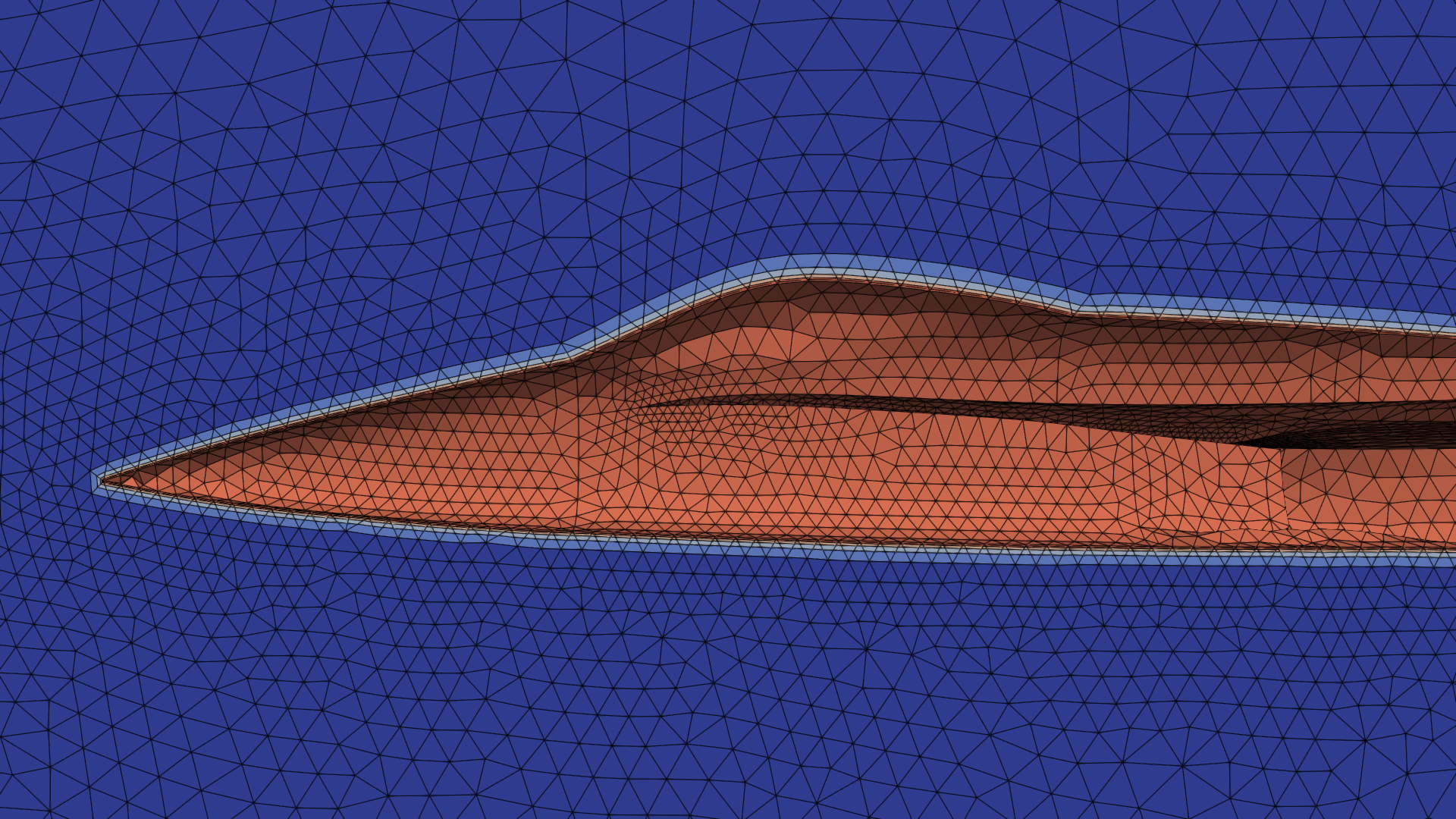}} \hfill
  \subfloat[Level 5]{\includegraphics[width=0.49\textwidth, trim={0, 5cm, 0, 5cm}, clip]{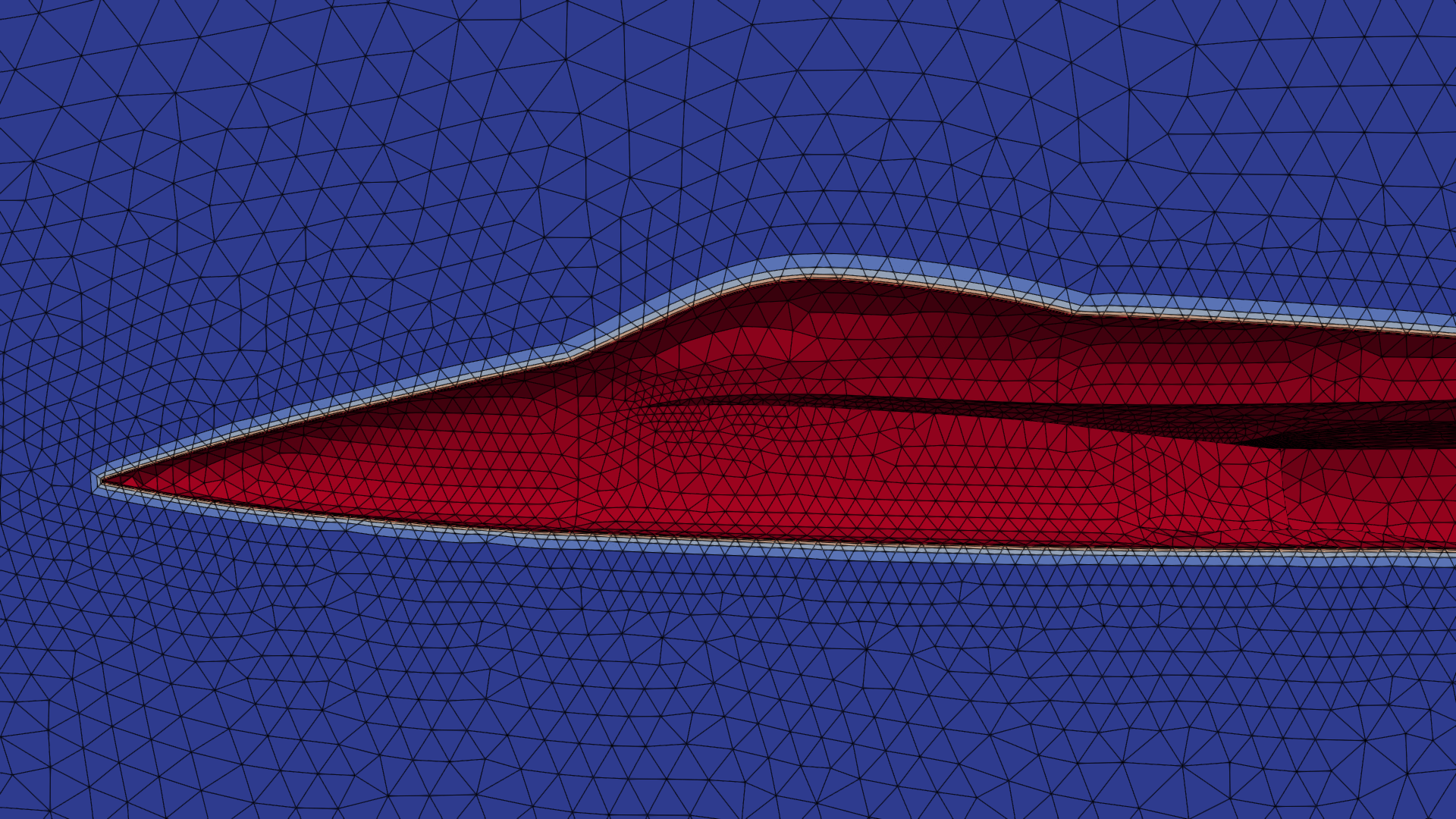}} \\
  \caption{YF-17 aircraft mesh. Lateral view of five levels of directional refinements.}
  \label{fig:YF17-lateral-refinements}
\end{figure}

\begin{figure}[!ht]
  \centering
  \includegraphics{f17-label-h.pdf}\\
  \subfloat[Level 1]{\includegraphics[width=0.49\textwidth, trim={0, 0, 0, 7cm}, clip]{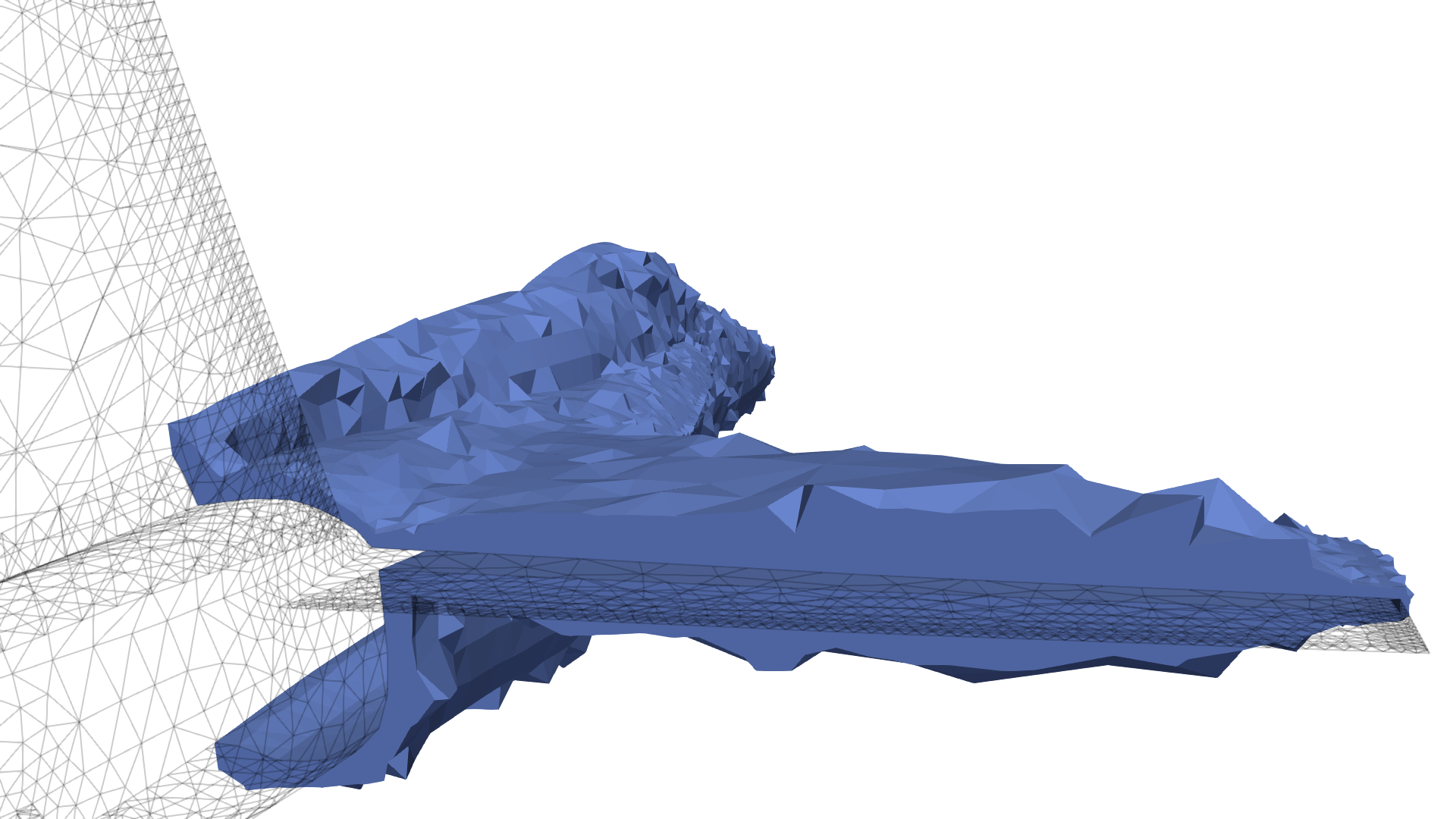}} \hfill
  \subfloat[Level 2]{\includegraphics[width=0.49\textwidth, trim={0, 0, 0, 7cm}, clip]{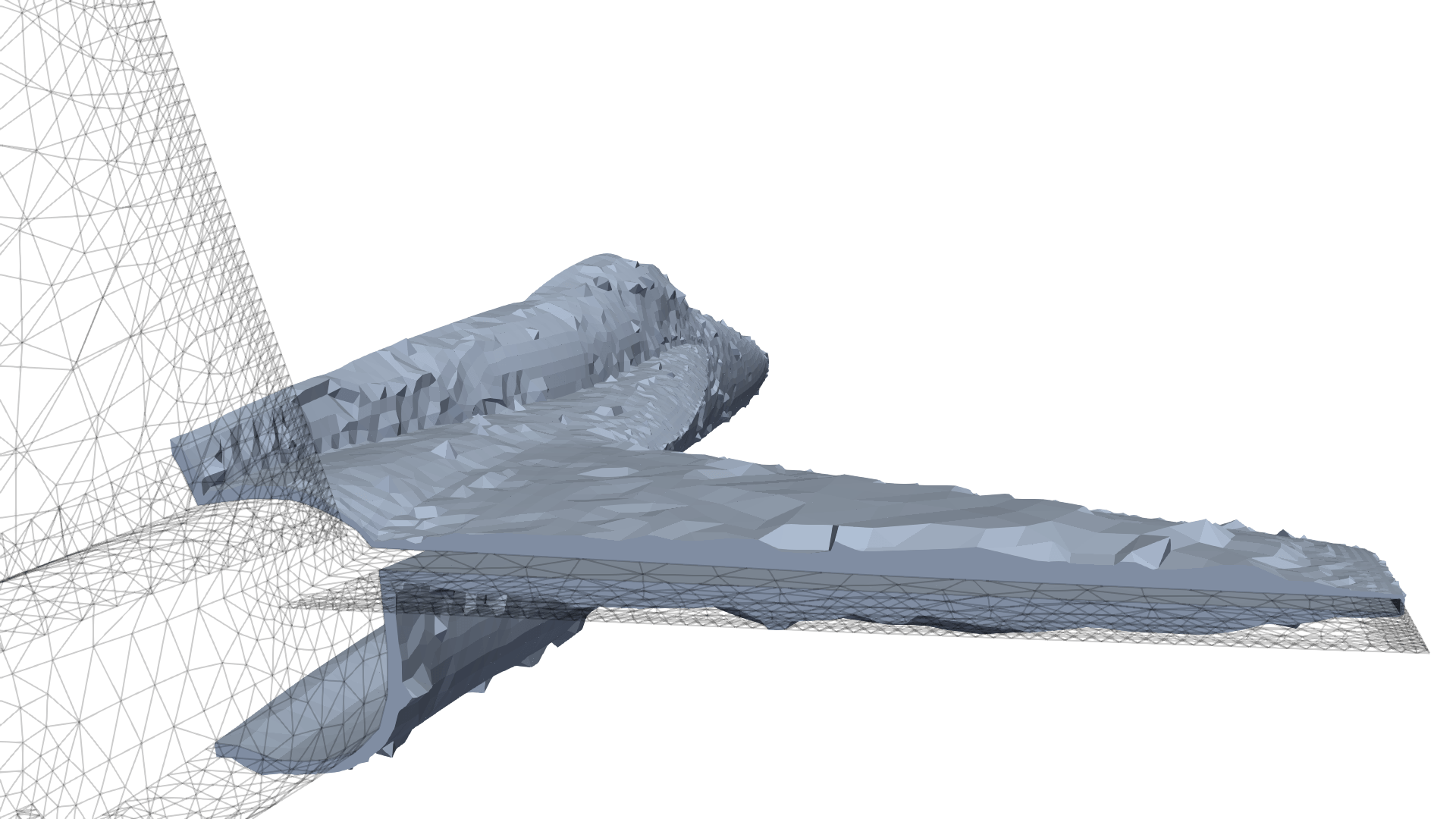}} \\
  \subfloat[Level 3]{\includegraphics[width=0.49\textwidth, trim={0, 0, 0, 7cm}, clip]{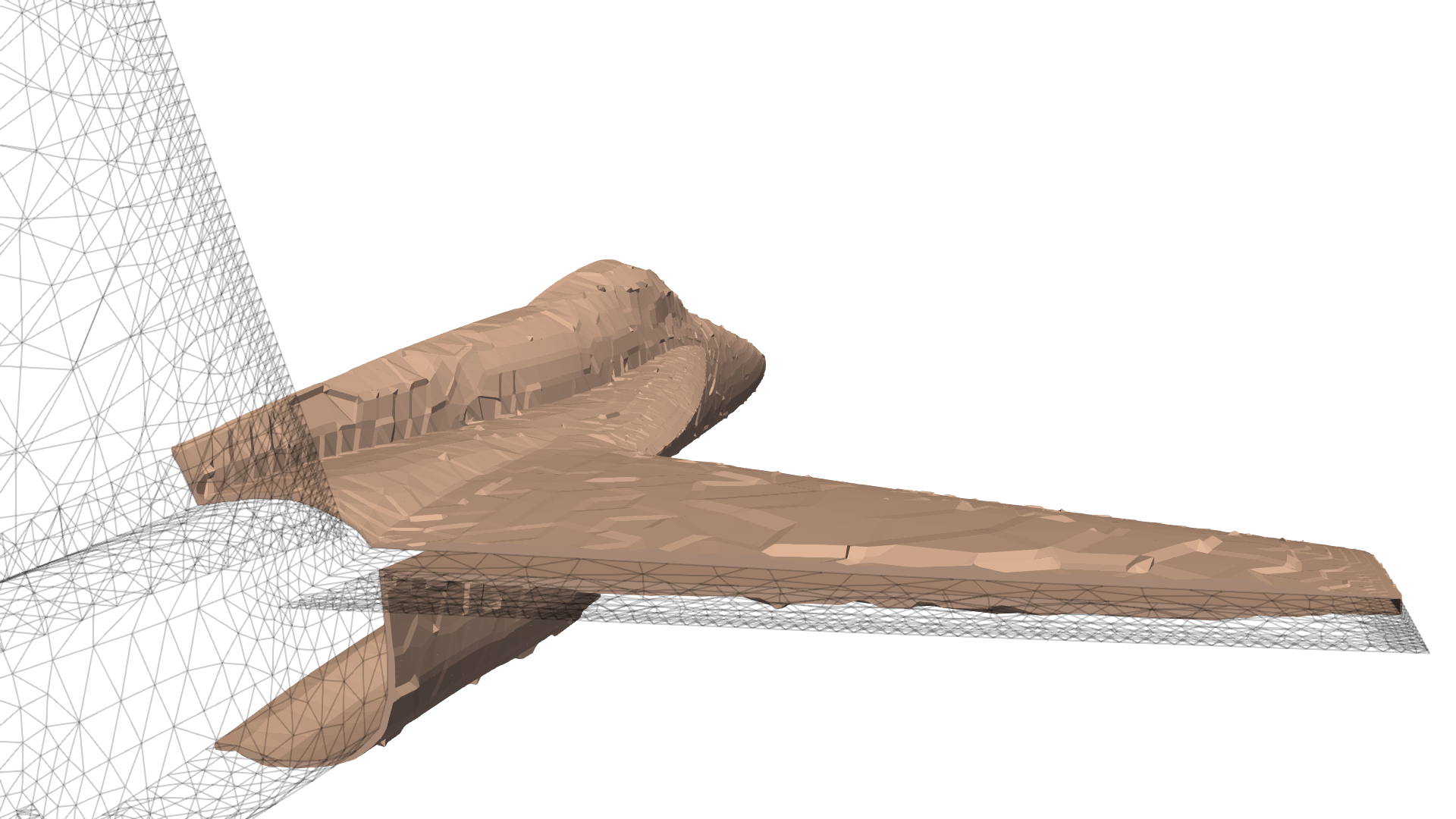}} \hfill
  \subfloat[Level 4]{\includegraphics[width=0.49\textwidth, trim={0, 0, 0, 7cm}, clip]{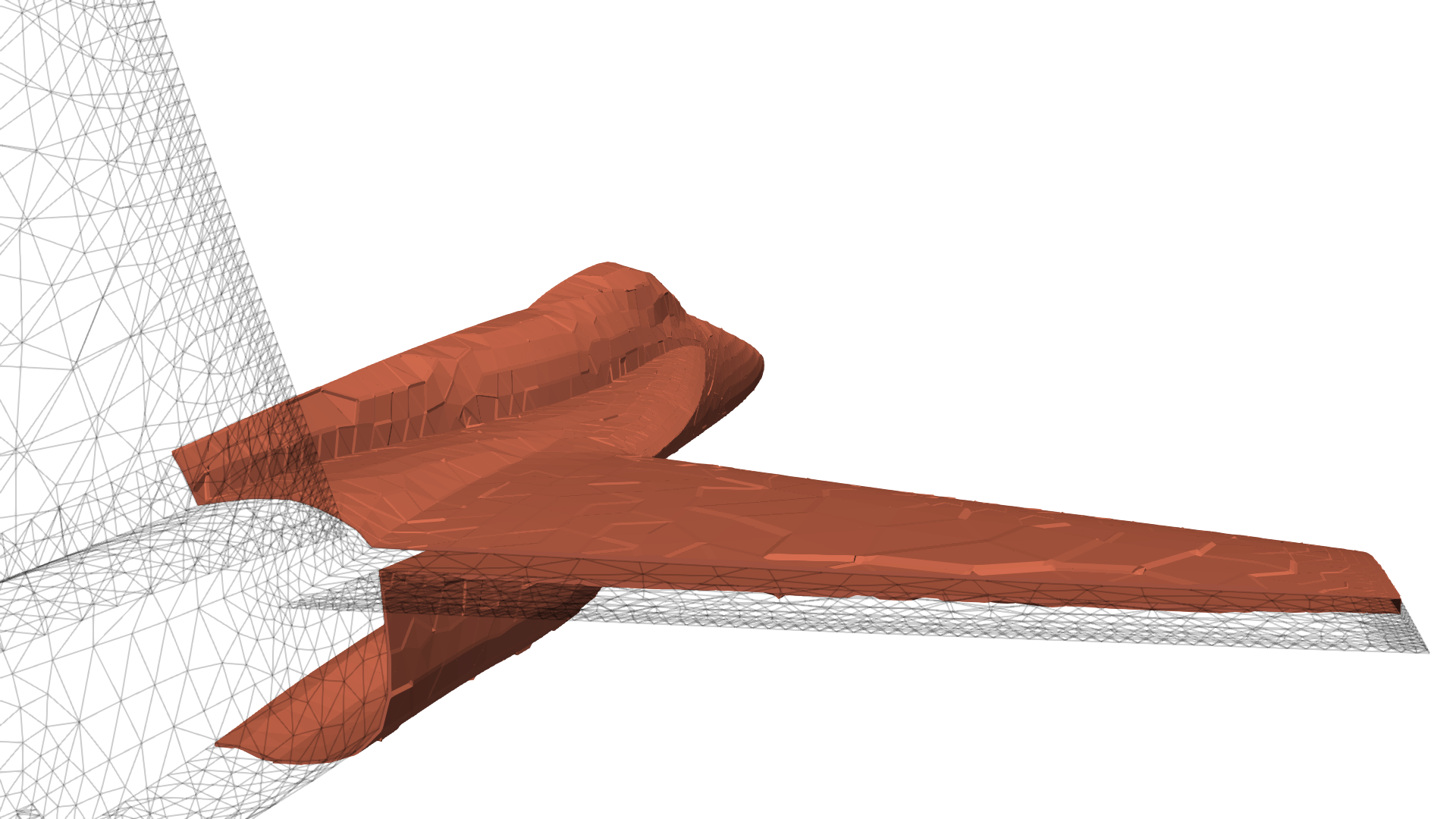}} \\
  \subfloat[Level 5]{\includegraphics[width=0.49\textwidth, trim={0, 0, 0, 7cm}, clip]{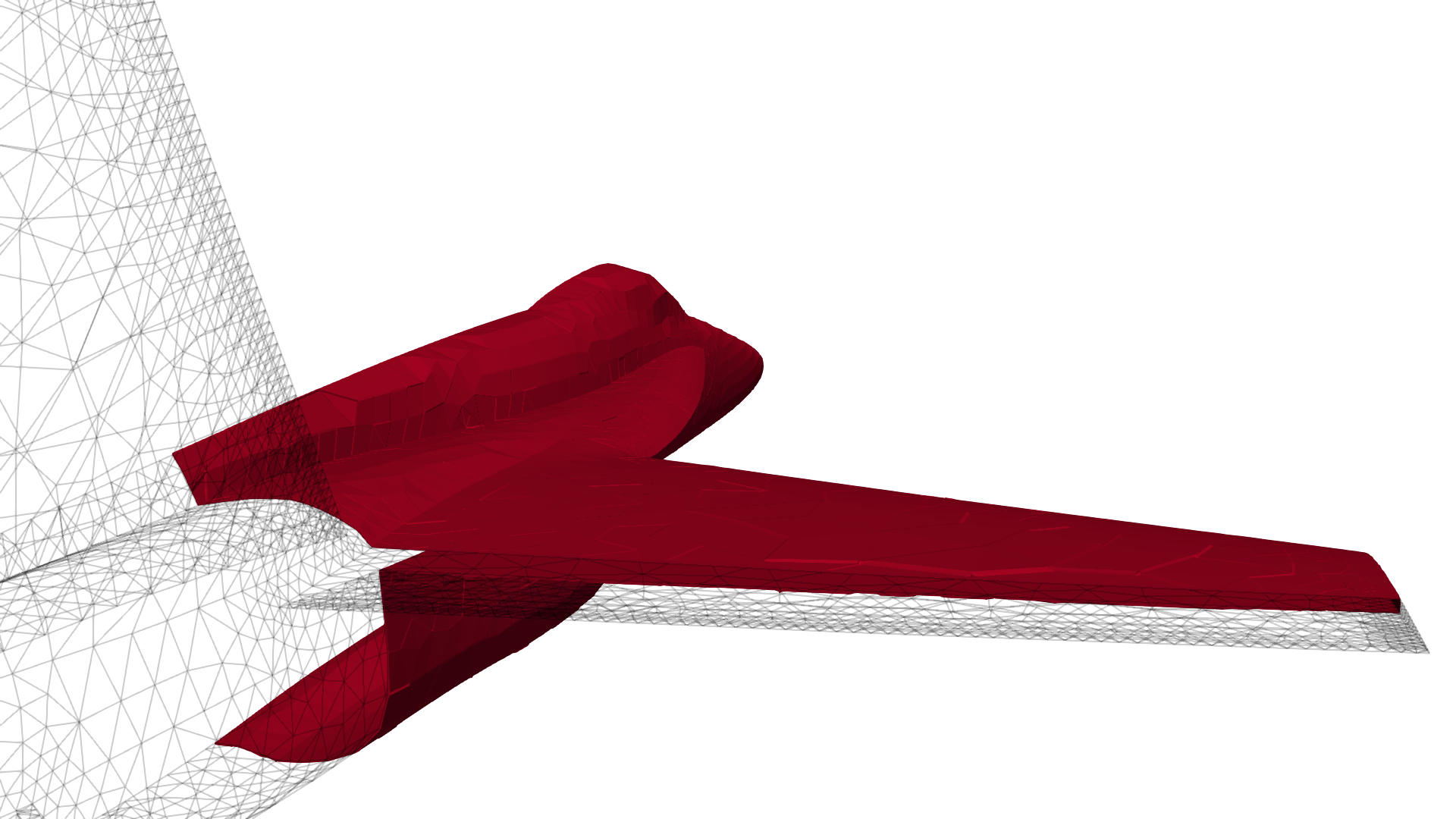}} \hfill
  \subfloat[All levels]{\includegraphics[width=0.49\textwidth, trim={0, 0, 0, 5cm}, clip]{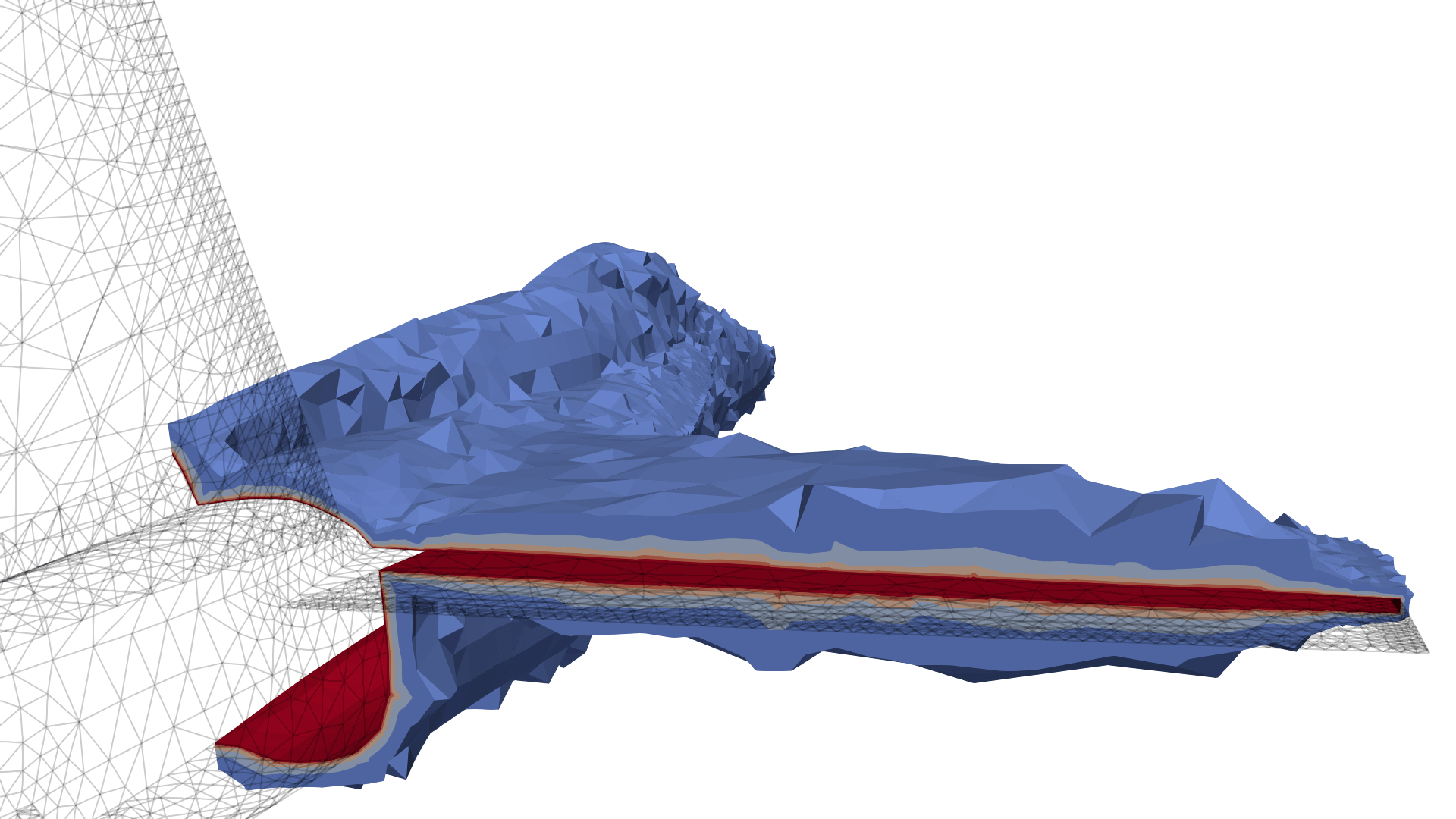}} \\
  \caption{YF-17 aircraft mesh. Rear view of five levels of directional refinements.}
  \label{fig:YF17-rear-refinements}
\end{figure}

In order to accurately capture the boundary layer effects of the airflow near the aircraft surface, a finer mesh is usually needed in that region. The directional refinement tool is tailored for this kind of problem, as it allows the user to refine elements that are neighbors to a given geometrical entity, such as a point, a curve or a surface. We apply five consecutive refinements, leading to the mesh layers displayed in Figures \ref{fig:YF17-lateral-refinements}-\ref{fig:YF17-rear-refinements}.

\begin{figure}[!ht]
  \centering
  \includegraphics[width=0.6\textwidth, valign=c]{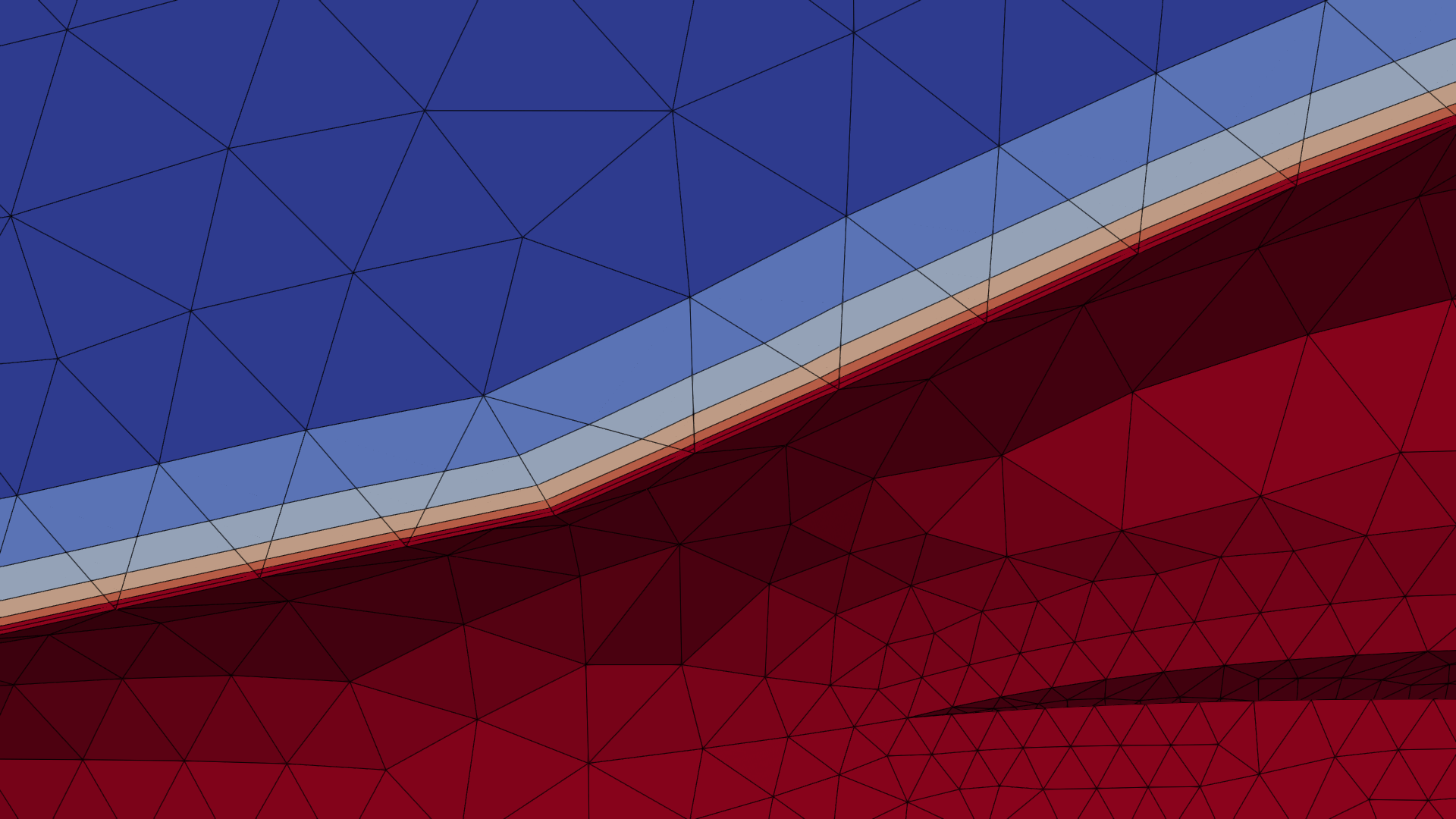} \quad
  \includegraphics[valign=c]{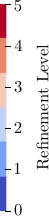}
  \caption{YF-17 aircraft mesh. Boundary layers after five levels of refinements.}
  \label{fig:YF17-lateral-zoom}
\end{figure}

\begin{figure}[!ht]
  \centering
  \subfloat[Level 0]{\includegraphics[width=0.3\textwidth, trim={0, 0, 0, 0cm}, clip, valign=c]{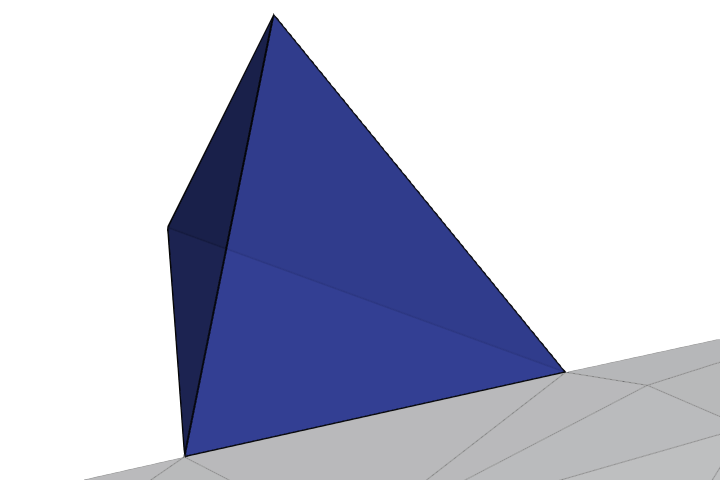}} \hfill
  \subfloat[Level 1]{\includegraphics[width=0.3\textwidth, trim={0, 0, 0, 0cm}, clip, valign=c]{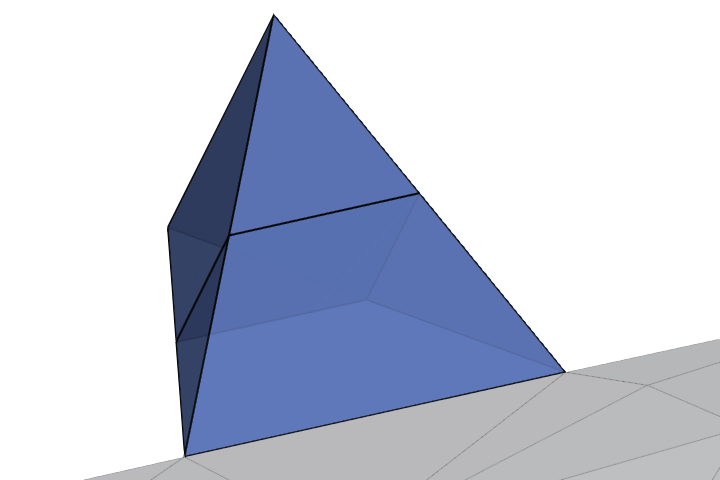}} \hfill
  \subfloat[Level 2]{\includegraphics[width=0.3\textwidth, trim={0, 0, 0, 0cm}, clip, valign=c]{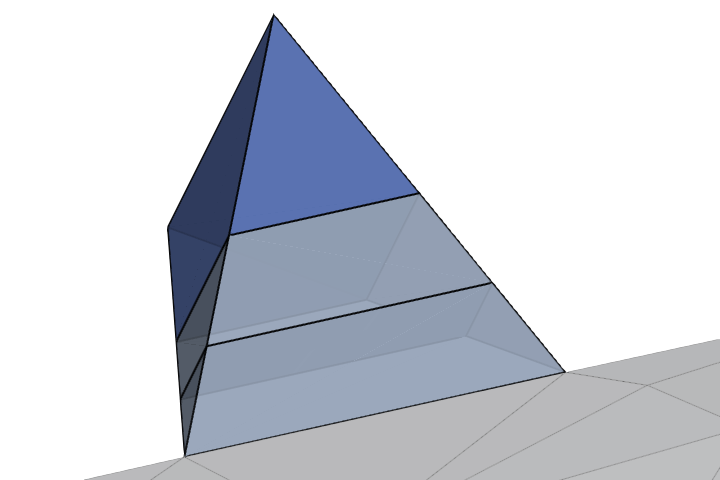}} \quad
  \includegraphics[valign=c]{f17-label-v.pdf} \\
  \subfloat[Level 3]{\includegraphics[width=0.3\textwidth, trim={0, 0, 0, 0cm}, clip, valign=c]{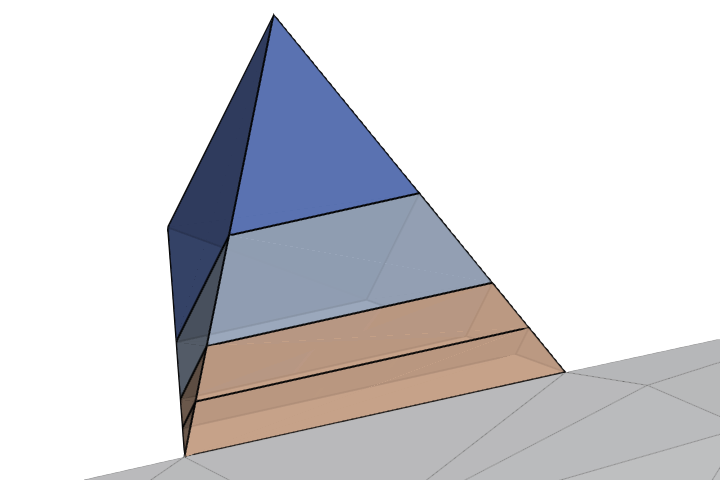}} \hfill
  \subfloat[Level 4]{\includegraphics[width=0.3\textwidth, trim={0, 0, 0, 0cm}, clip, valign=c]{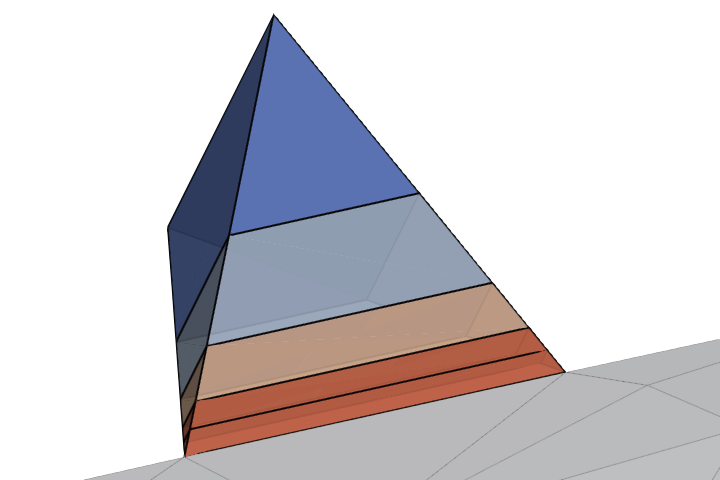}} \hfill
  \subfloat[Level 5]{\includegraphics[width=0.3\textwidth, trim={0, 0, 0, 0cm}, clip, valign=c]{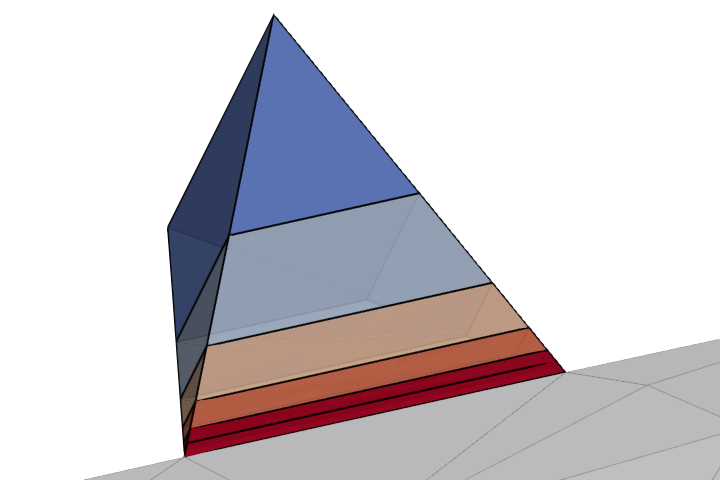}} \quad
  \includegraphics[valign=c]{f17-label-v.pdf} \\
  \caption{YF-17 aircraft mesh. Successive refinements over an element that shares a line with the aircraft surface.}
  \label{fig:YF17-touching-line}
\end{figure}

\begin{figure}[!ht]
  \centering
  \subfloat[Level 0]{\includegraphics[width=0.3\textwidth, trim={0, 0, 0, 0cm}, clip, valign=c]{f17-el-surface-0.png}} \hfill
  \subfloat[Level 1]{\includegraphics[width=0.3\textwidth, trim={0, 0, 0, 0cm}, clip, valign=c]{f17-el-surface-1.png}} \hfill
  \subfloat[Level 2]{\includegraphics[width=0.3\textwidth, trim={0, 0, 0, 0cm}, clip, valign=c]{f17-el-surface-2.png}} \quad
  \includegraphics[valign=c]{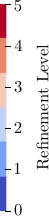} \\
  \subfloat[Level 3]{\includegraphics[width=0.3\textwidth, trim={0, 0, 0, 0cm}, clip, valign=c]{f17-el-surface-3.png}} \hfill
  \subfloat[Level 4]{\includegraphics[width=0.3\textwidth, trim={0, 0, 0, 0cm}, clip, valign=c]{f17-el-surface-4.png}} \hfill
  \subfloat[Level 5]{\includegraphics[width=0.3\textwidth, trim={0, 0, 0, 0cm}, clip, valign=c]{f17-el-surface-5.png}} \quad
  \includegraphics[valign=c]{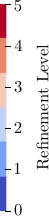} \\
  \caption{YF-17 aircraft mesh. Successive refinements over an element that shares a facet with the aircraft surface.}
  \label{fig:YF17-touching-surface}
\end{figure}

In Figure \ref{fig:YF17-lateral-zoom}, it can be observed that the aircraft surface mesh is not modified as new layers of elements are introduced in the 3D domain surrounding the surface. As stated in Section \ref{sec:def-ref-pat}, the refinement at the element level occurs following a predefined pattern stored in a database. For tetrahedra, the pattern that matches the desired refinement strategy will not always produce child elements whose topology is the same as the parent element. To better exemplify this, figures \ref{fig:YF17-touching-line}-\ref{fig:YF17-touching-surface} depict two tetrahedral elements that share a line and a facet, respectively, with the aircraft surface. In the first case, the first level of refinement leads to two child elements with wedge topology, while in the second case, the tetrahedron is split into a smaller tetrahedron and a wedge. During the third and fourth refinement steps, some hexahedral elements are created when the father element is refined towards a line, while for the facet case, the remaining refinements produce only wedge elements.

\subsection{Fracture propagation}

\begin{figure}[!ht]
  \centering
  \subfloat[]{%
    \includegraphics[width=0.45\textwidth, valign=c]{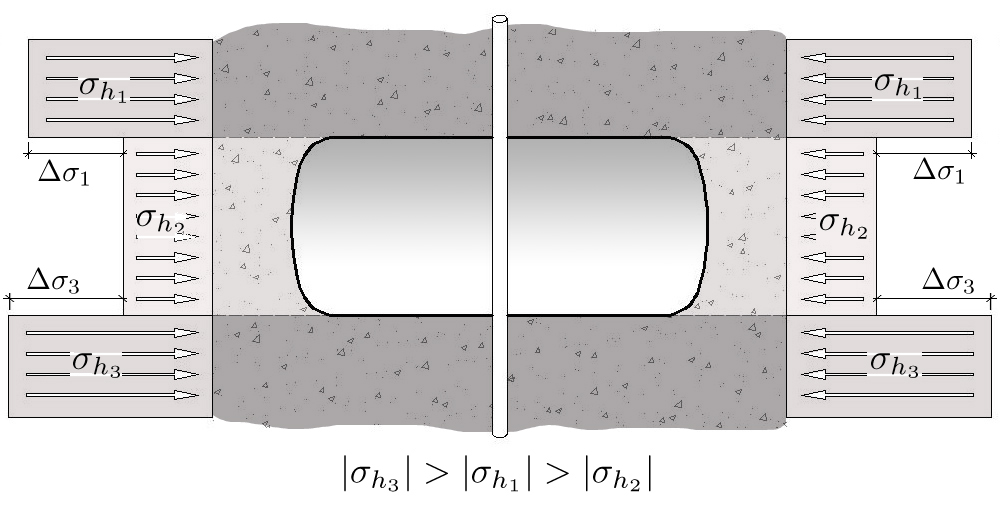}%
  } \qquad
  \subfloat[]{%
    \includegraphics[width=0.4\textwidth, valign=c]{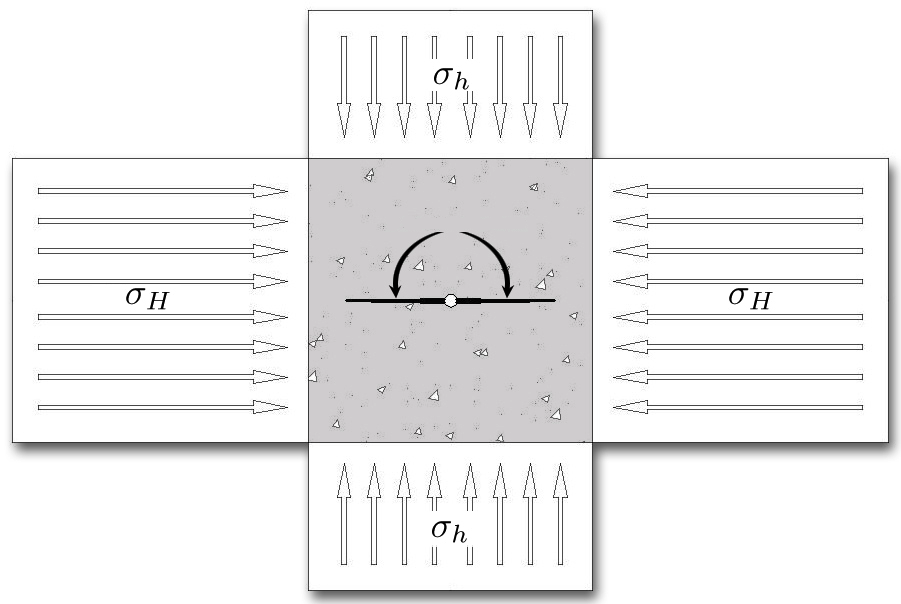}%
  }
  \caption{Illustration of hydraulic fracture propagation. (a) Side view and (b) top view with arrows indicating the 2 fractures being propagated from the wellbore.}
  \label{fig:hf-problem}
\end{figure}

\begin{figure}[!ht]
  \centering
  \subfloat[]{%
    \includegraphics[width=0.45\textwidth, valign=c]{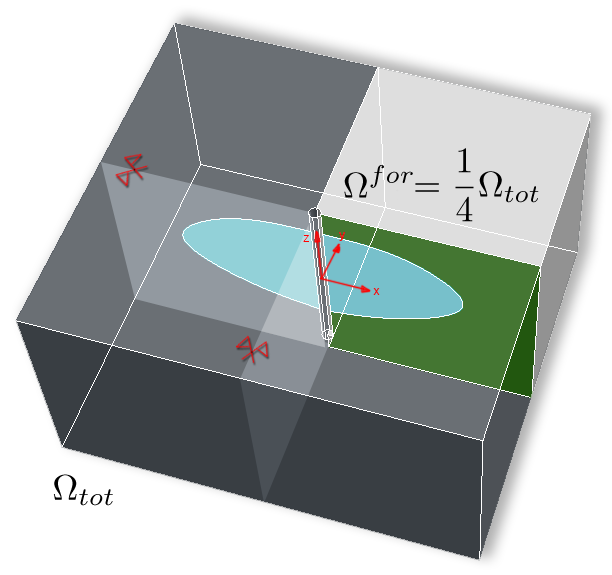}%
  } \qquad
  \subfloat[]{%
    \includegraphics[width=0.3\textwidth, valign=c]{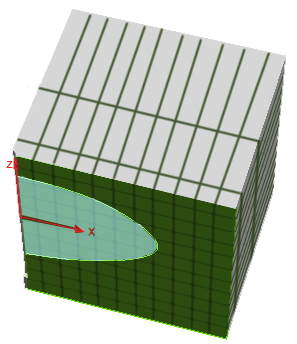}%
  }
  \caption{Illustration of the simulation domain used within the context of the FEM. (a) Complete domain and (b) one-fourth of the domain used for the simulation.}
  \label{fig:hf-domain}
\end{figure}

Hydraulic fracturing is a common procedure used in the oil industry to enhance the rates of injection and production of a wellbore. In this process, a fracture is propagated through the injection of pressurized fluid through the wellbore. This procedure is illustrated for the case of a vertical wellbore in Figure \ref{fig:hf-problem}, where two vertical fractures are propagated perpendicular to the minimum in-situ stress direction. Several methodologies have been proposed to simulate this process within the context of the Finite Element Method \cite{2018_Nathan_HF_Alg_Improv,2019_Nathan_validation,2021_Nathan_multiple_HF,DEVLOO2006,Adachi2007,Lecampion_hydraulic_fracture_review_2018,MearMultipleCrackGrowth,SALIMZADEH20179,WU2017,WaismanHF2019}. One possibility to increase efficiency is to only model one-fourth of the domain as shown in Figure \ref{fig:hf-domain}.

\begin{figure}[!ht]
  \centering
  \subfloat[]{%
    \includegraphics[width=0.35\textwidth, valign=c]{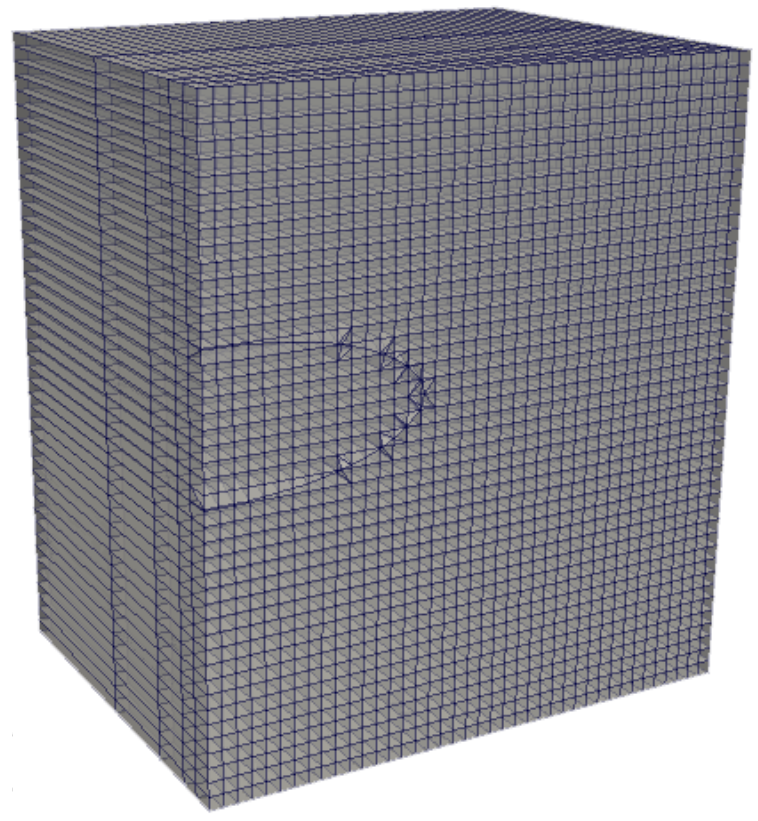}%
  } \qquad\qquad
  \subfloat[]{%
    \includegraphics[width=0.4\textwidth, valign=c]{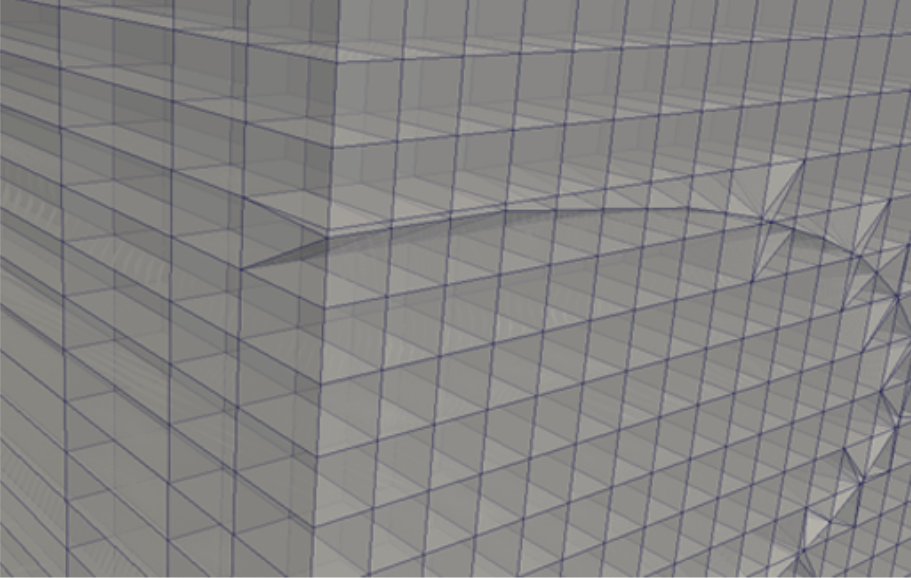}%
    \vphantom{\includegraphics[width=0.4\textwidth,valign=c]{porous-media-initial-mesh.png}}
  }
  \caption{Rock matrix with a vertical fracture. (a) Initial mesh and (b) fracture detail.}
  \label{fig:porous-media-initial-mesh}
\end{figure}

The basis of this problem is to define a line that represents the geometry of the fracture. This line intersects the 3D elements of the domain, and the refinement patterns are used to create a mesh conforming to this line as shown in Figure \ref{fig:porous-media-initial-mesh}.

\begin{figure}[!ht]
  \centering
  \includegraphics[width=0.8\textwidth]{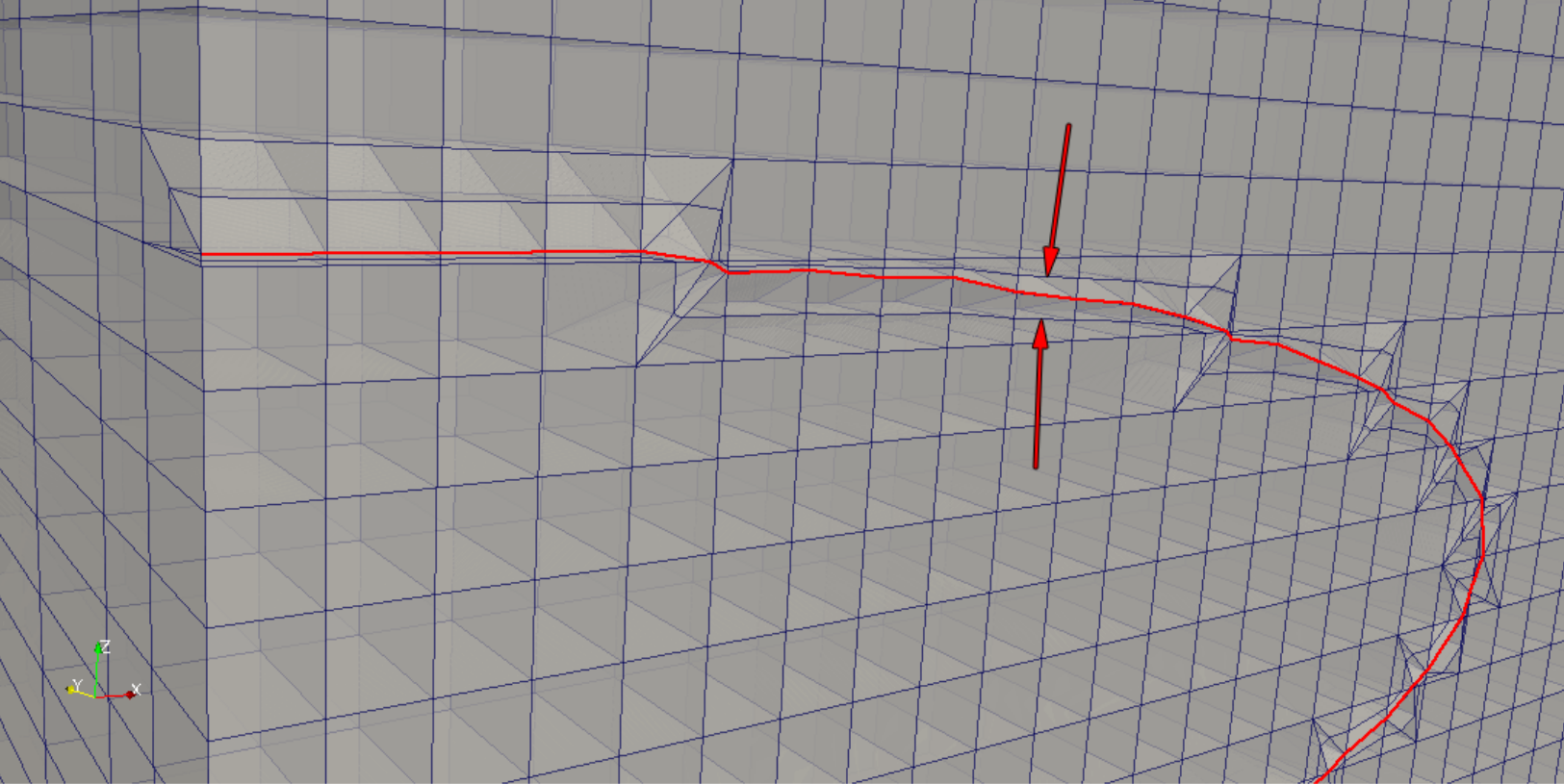}
  \caption{Rock matrix with a vertical fracture. Refinement towards the fracture contour.}
  \label{fig:porous-media-refined-mesh}
\end{figure}

It is well known that this kind of problem presents a singular stress field in the vicinity of the crack contour, so one possible strategy to tackle this issue is to refine the mesh near this region in order to better capture this stress concentration. To avoid global mesh refinement and improve the computational efficiency, one directional mesh refinement step is applied towards the fracture contour, as shown in Figure \ref{fig:porous-media-refined-mesh}.

\subsection{A structural casket designed to run experiments of petroleum engineering}

\begin{figure}[!ht]
  \centering
  \subfloat[]{%
    \includegraphics[width=0.35\textwidth, valign=c]{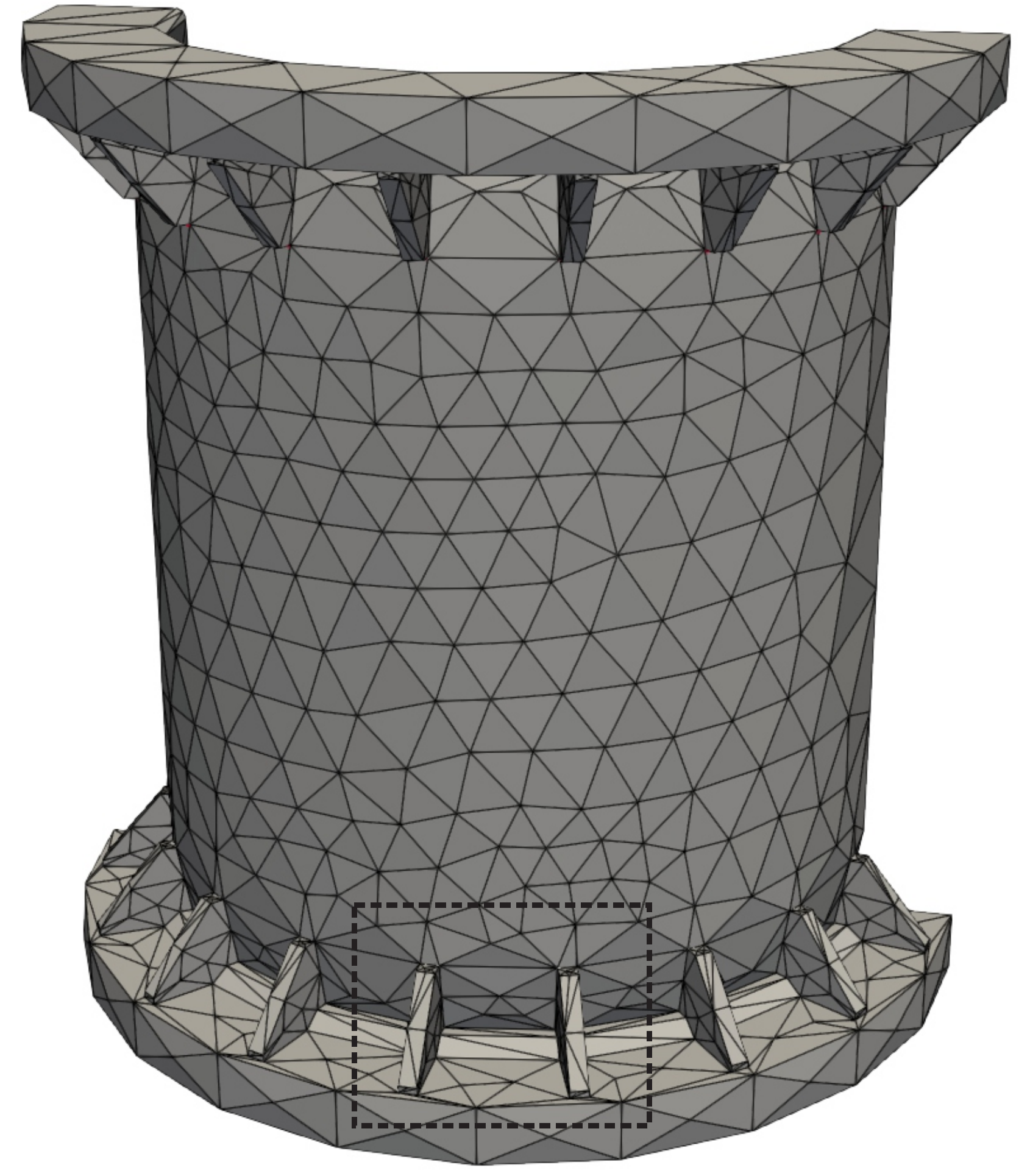}%
  } \qquad
  \subfloat[\label{fig:casket-stiffners}]{%
    \includegraphics[width=0.35\textwidth, valign=c]{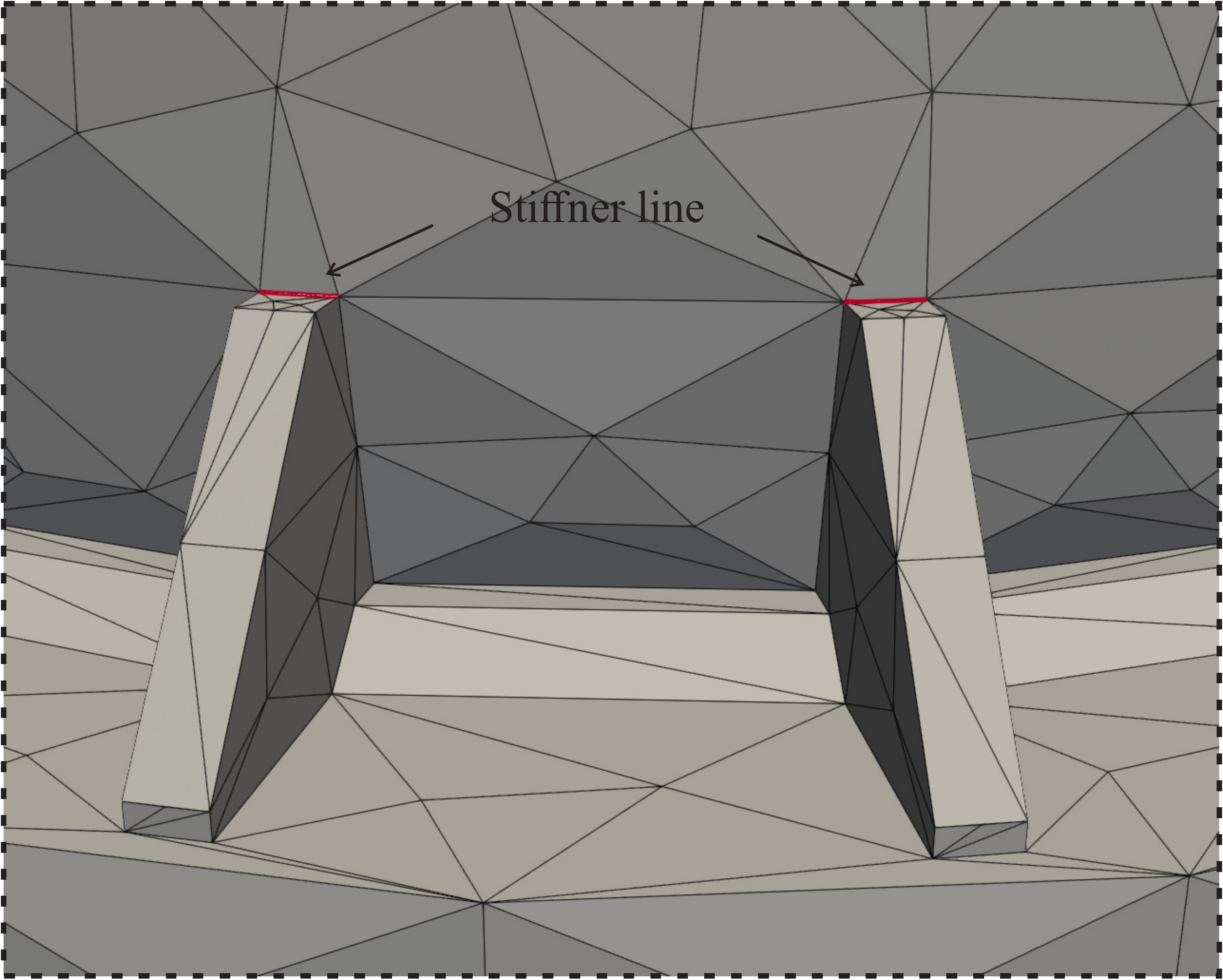}%
    \vphantom{\includegraphics[width=0.35\textwidth,valign=c]{perspective-view.pdf}}
  }
  \caption{Structural casket. (a) Initial mesh and (b) detail of the stiffeners.}
  \label{fig:casket-initial-mesh}
\end{figure}

This example is based on a real prototype being developed at Unicamp. It consists of a structural module made of steel, designed to withstand a uniform pressure load on its inner surface. The geometry and boundary conditions lead to an axisymmetric problem, so for simplicity, only half of the casket is modeled according to Figure \ref{fig:casket-initial-mesh}. Stiffners are placed on the upper and lower flanges to prevent buckling and increase local stability. However, the presence of these stiffeners induces a stress concentration that needs to be accurately captured by the numerical simulation. When using FEM formulations, a major concern is the mesh quality in this region which generally needs to be more refined than the rest of the domain.

\begin{figure}[!ht]
  \centering
  \includegraphics{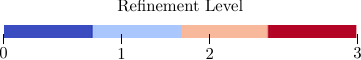} \\
  \subfloat[Level 0]{%
    \includegraphics[width=0.24\textwidth, valign=c]{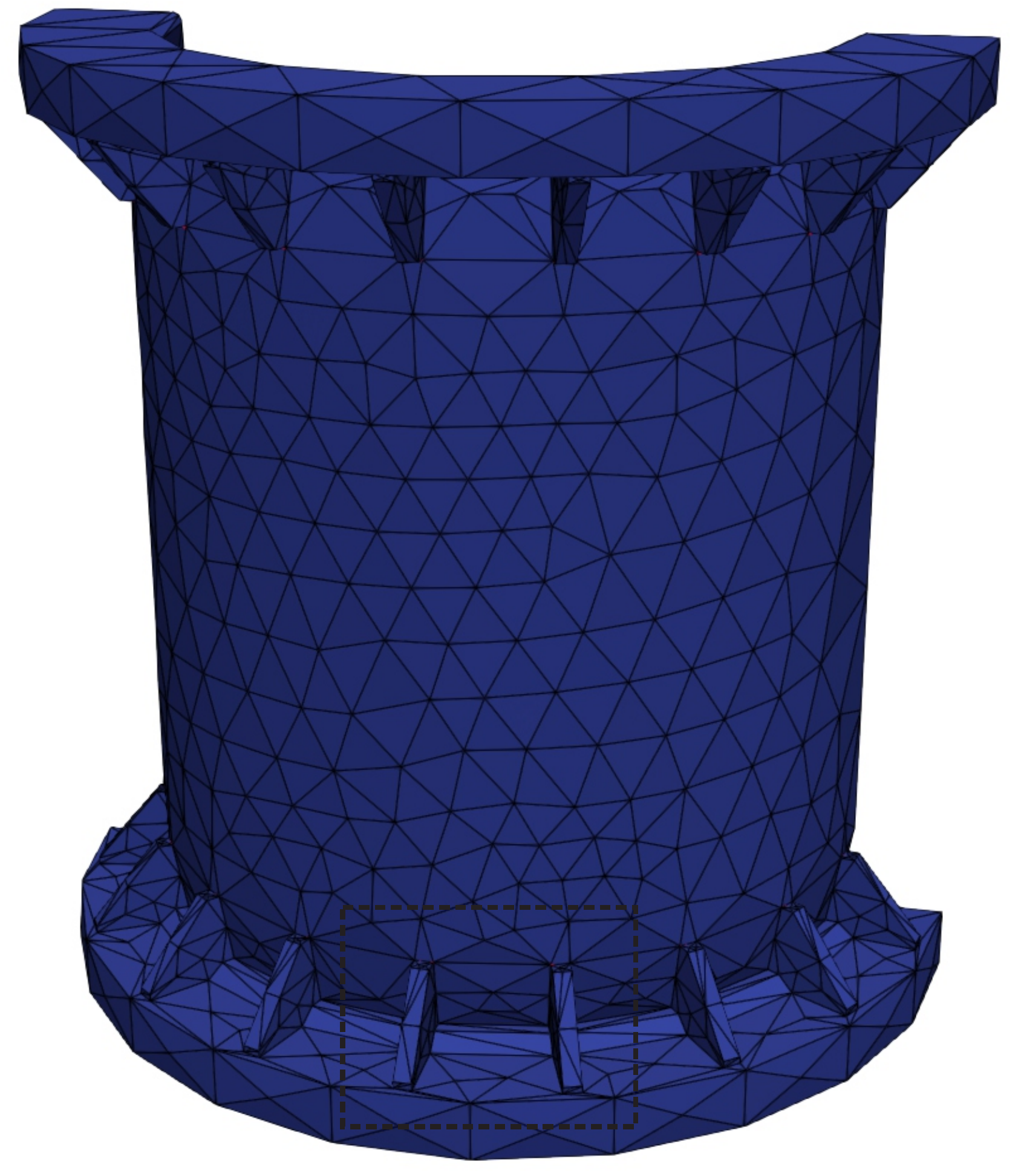}%
    \includegraphics[width=0.24\textwidth, valign=c]{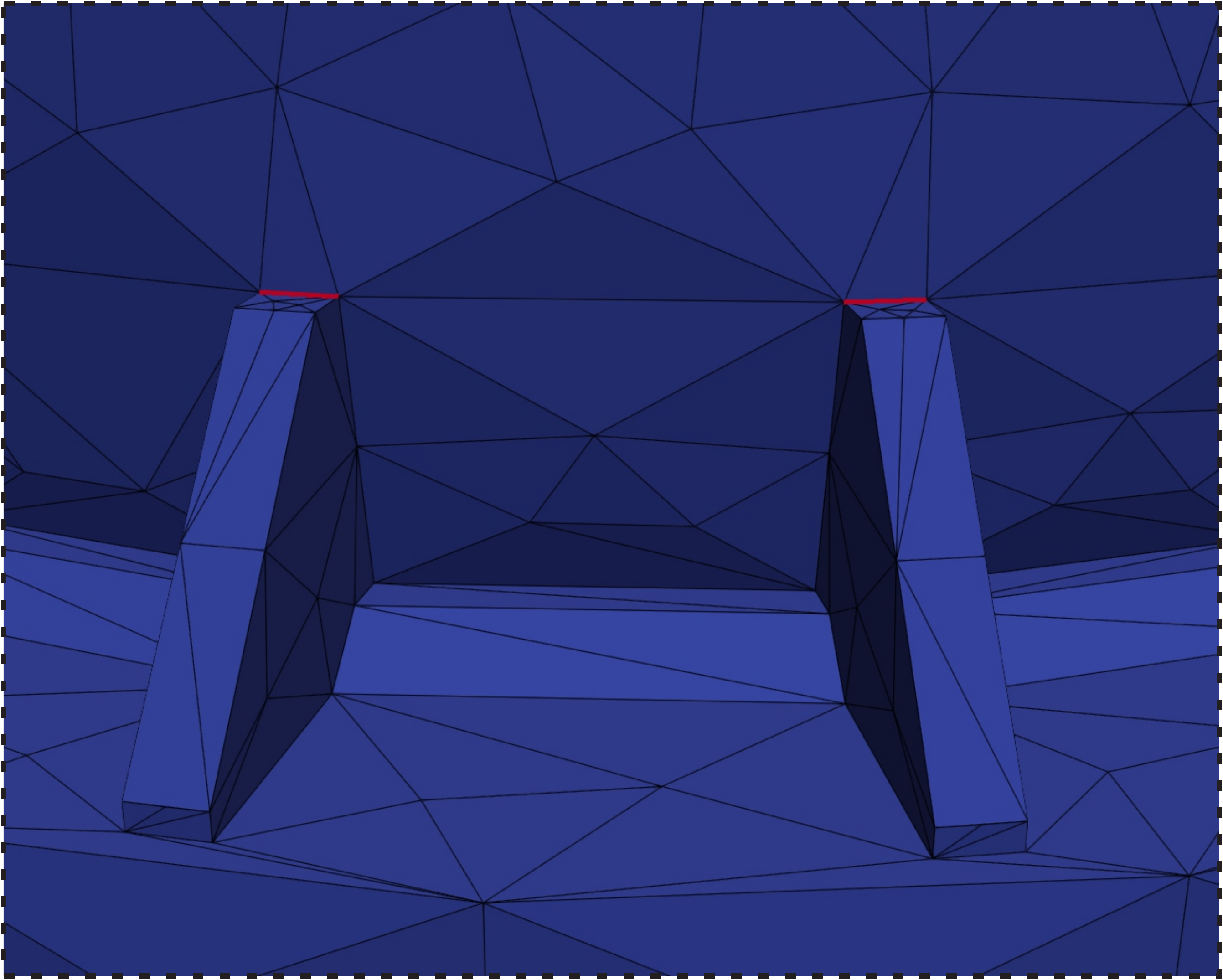}%
  } \hfill
  \subfloat[Level 1]{%
    \includegraphics[width=0.24\textwidth, valign=c]{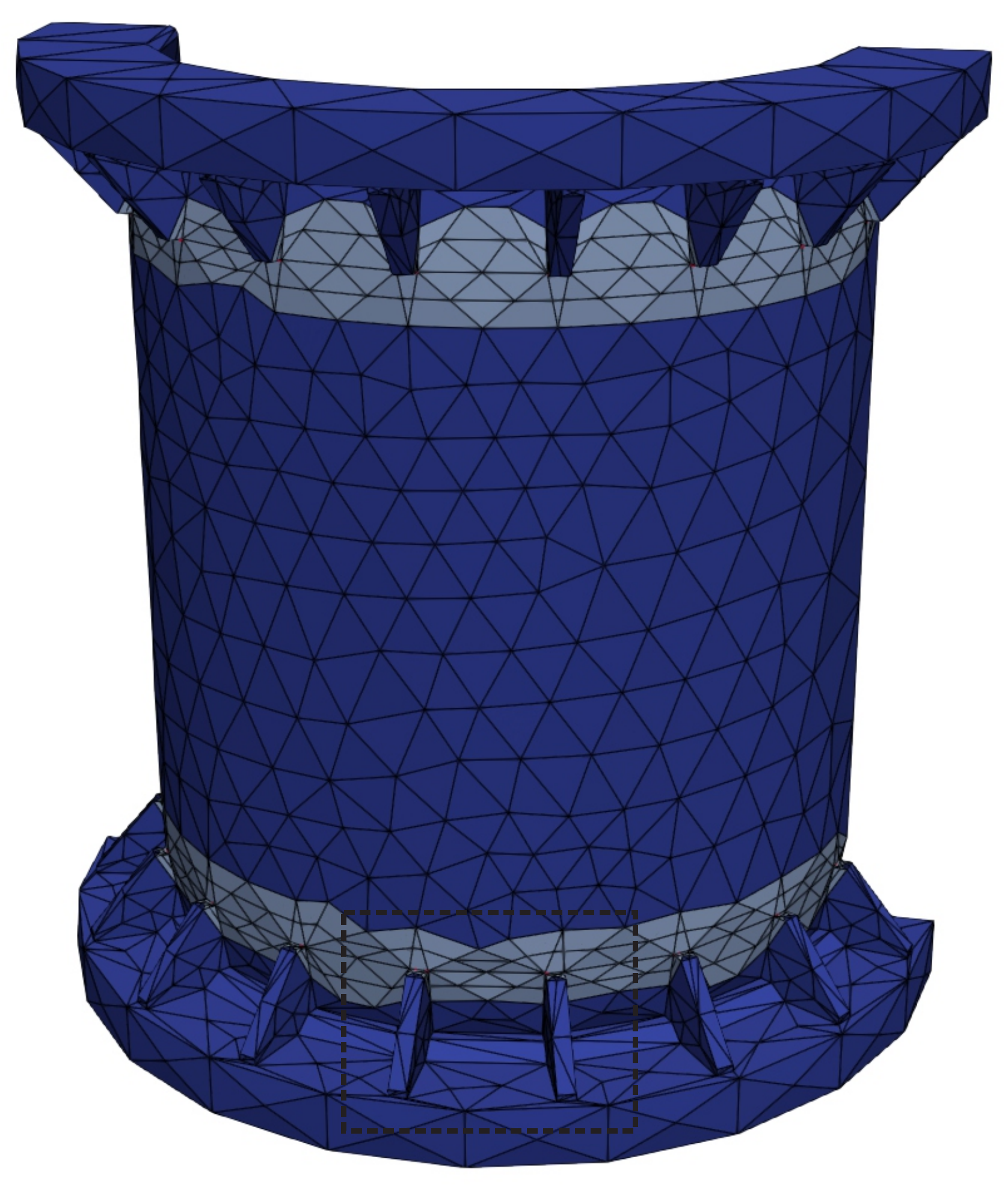}%
    \includegraphics[width=0.24\textwidth, valign=c]{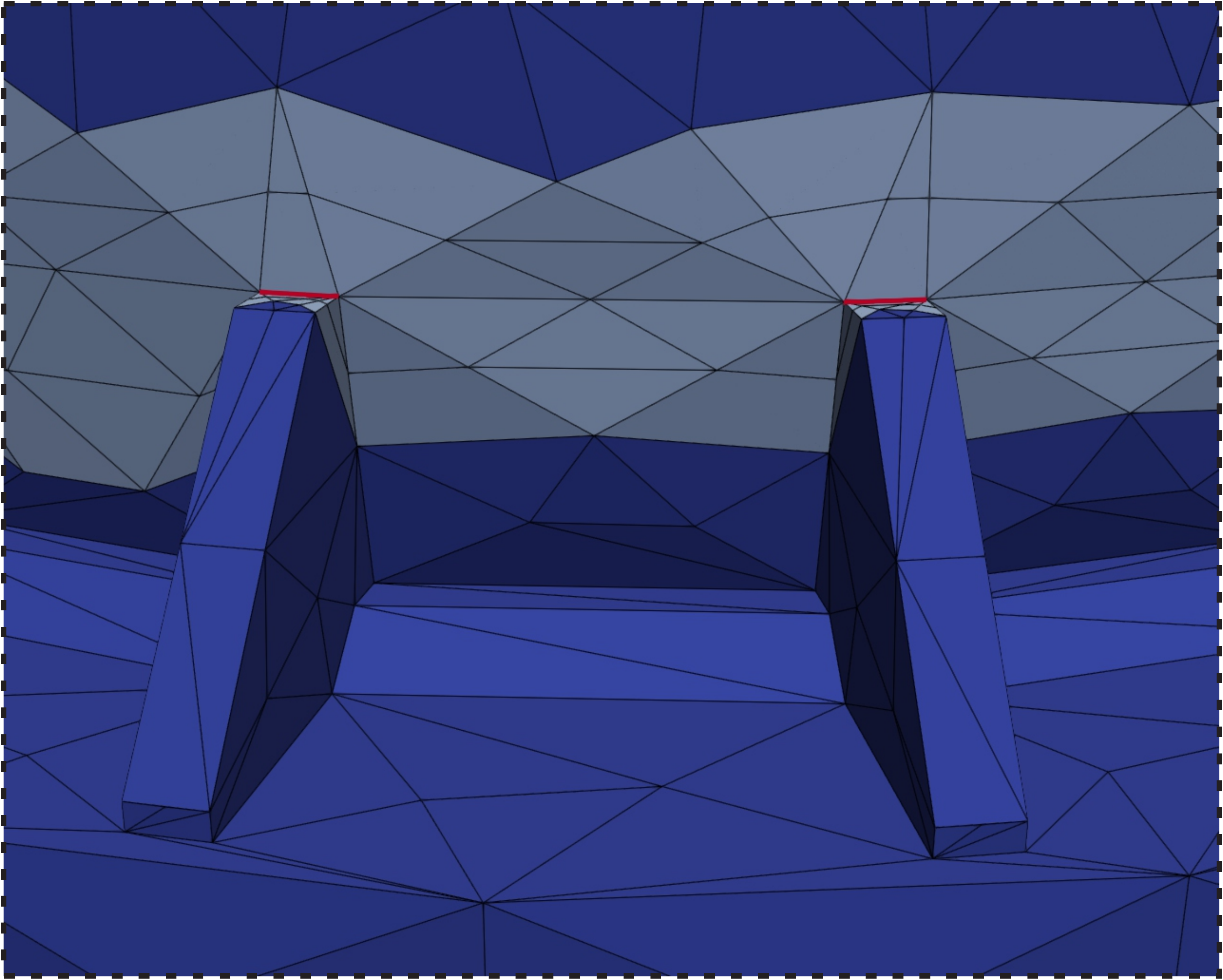}%
  } \\
  \subfloat[Level 2]{%
    \includegraphics[width=0.24\textwidth, valign=c]{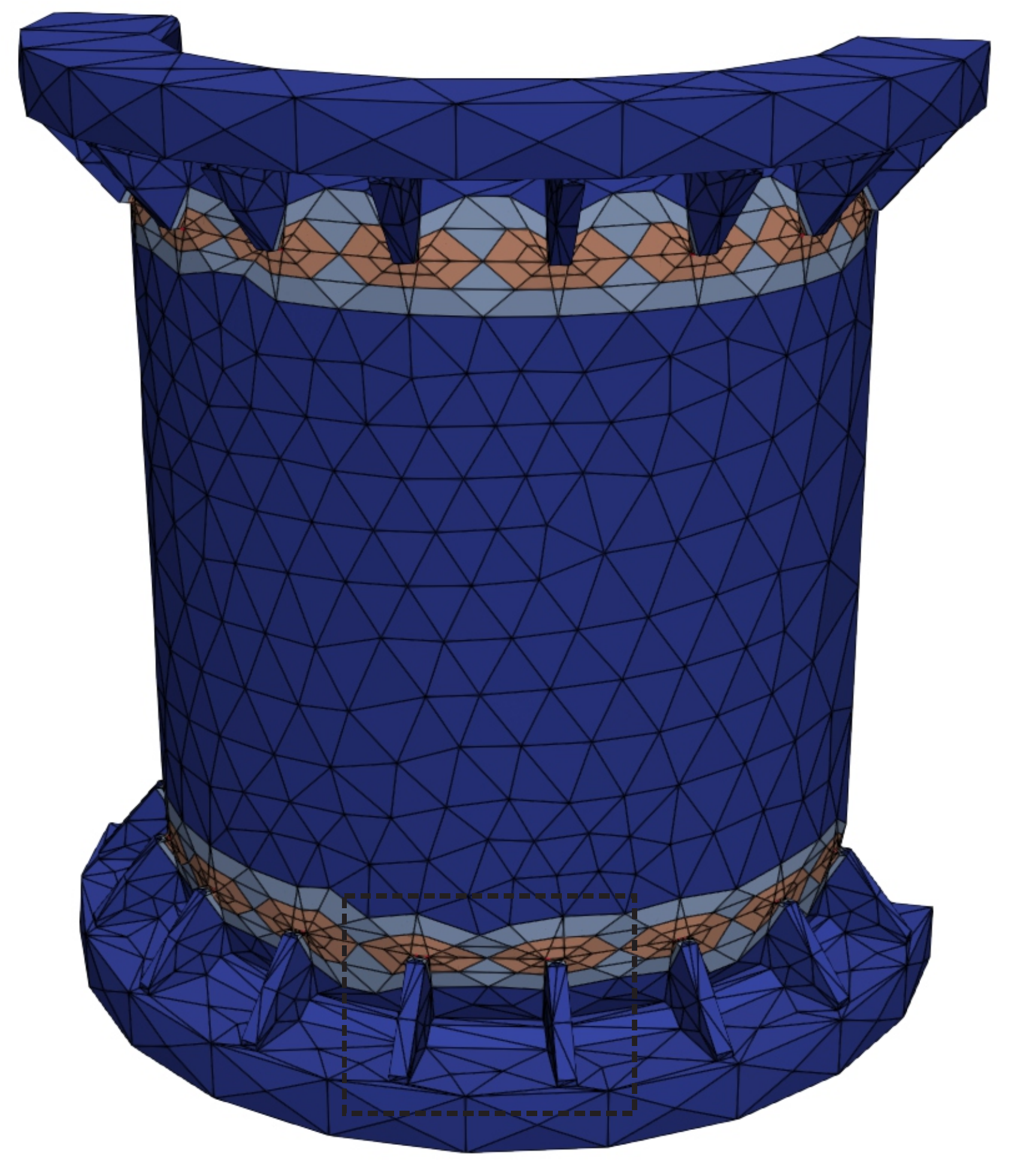}%
    \includegraphics[width=0.24\textwidth, valign=c]{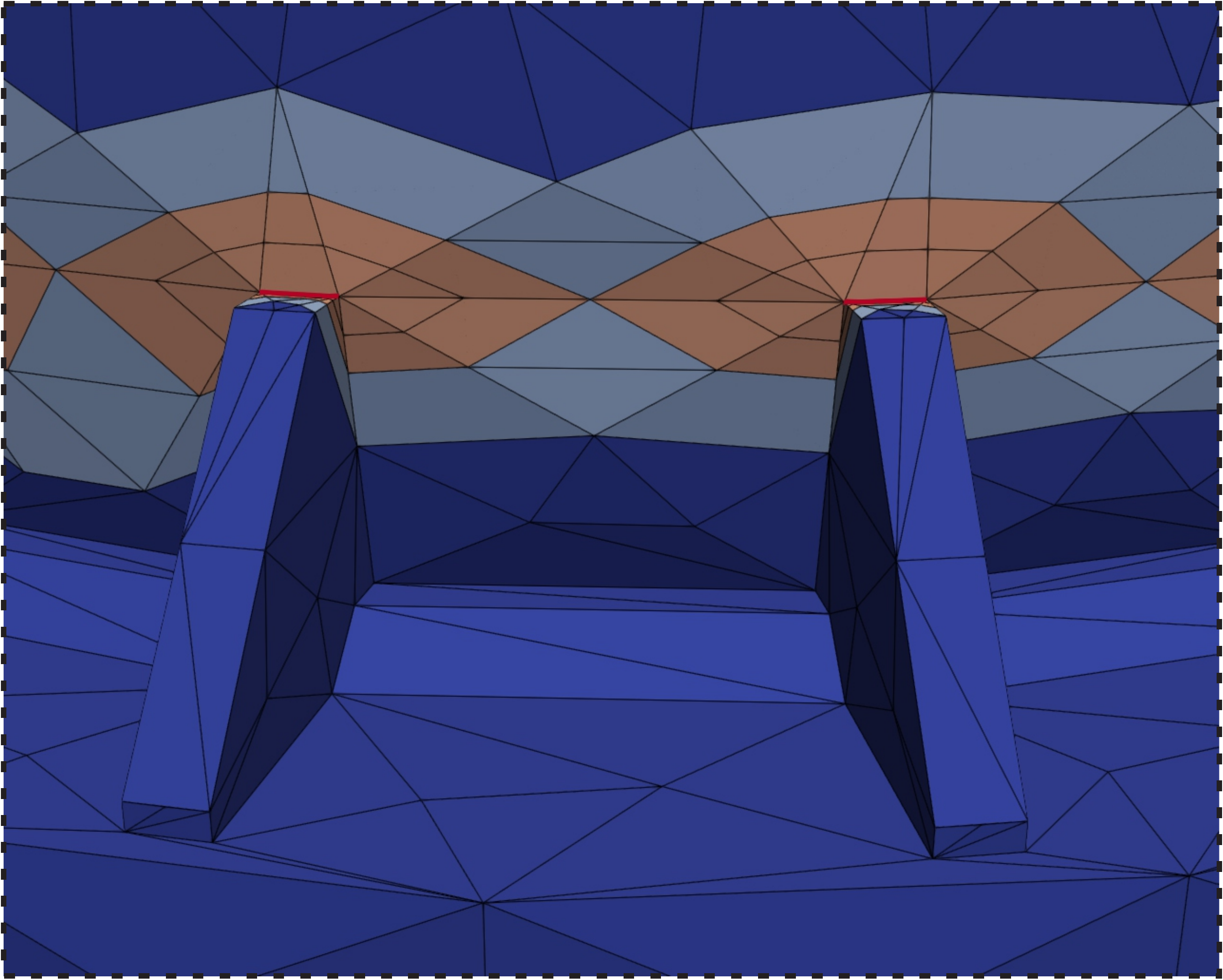}%
  } \hfill
  \subfloat[Level 3]{%
    \includegraphics[width=0.24\textwidth, valign=c]{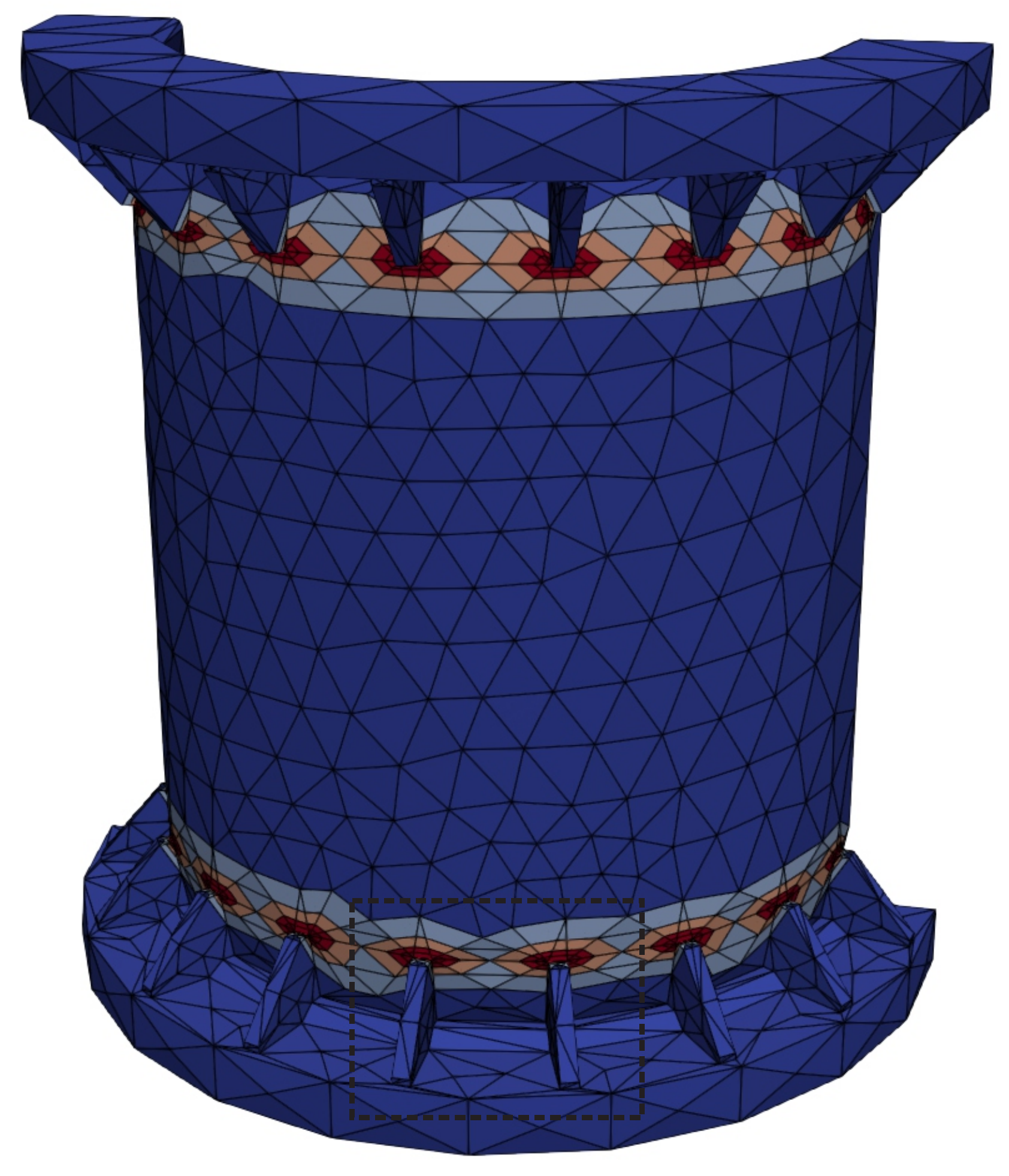}%
    \includegraphics[width=0.24\textwidth, valign=c]{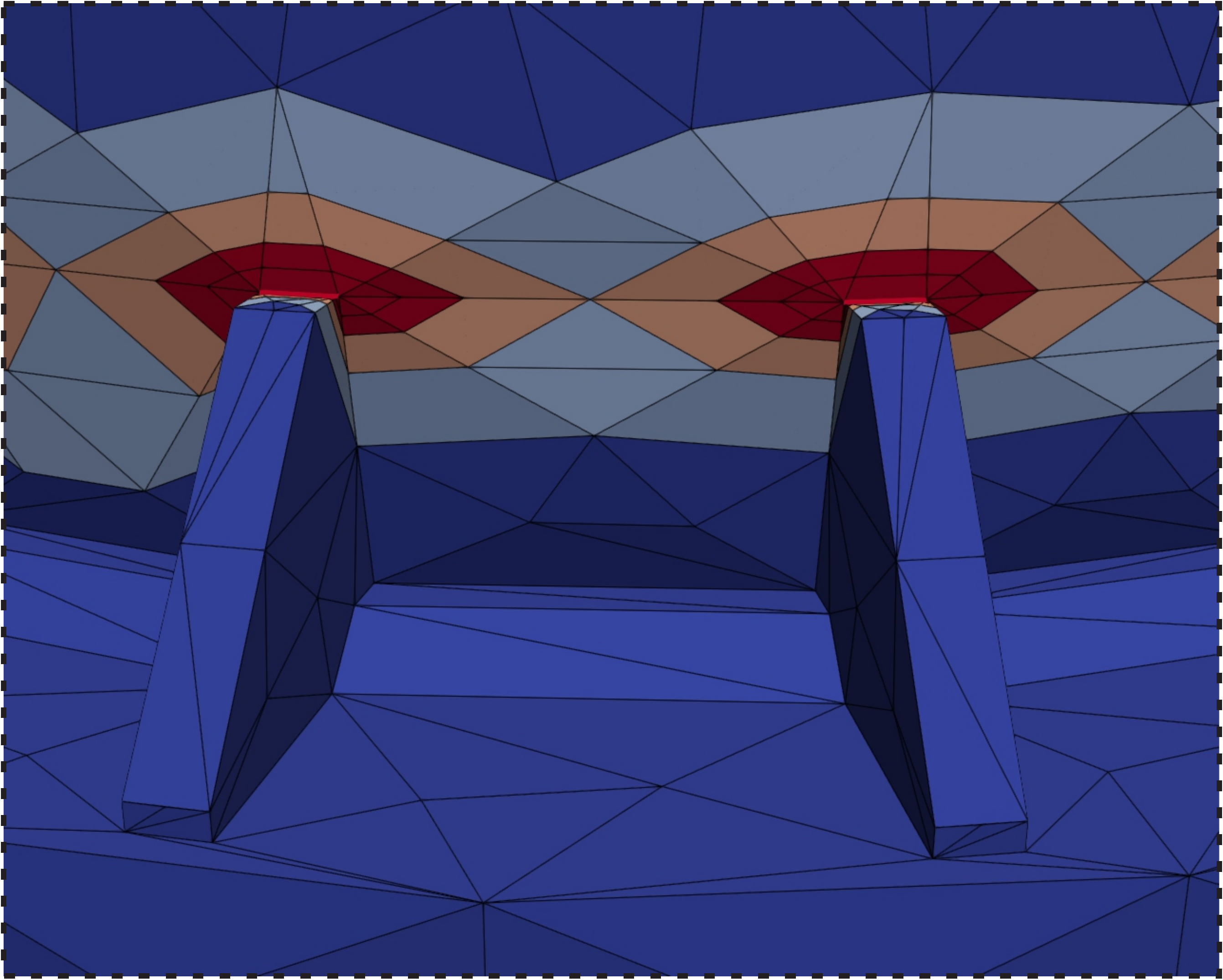}%
  }
  \caption{Structural casket. Successive refinements towards the stiffners}
  \label{fig:casket-refined-mesh}
\end{figure}

One option to tackle this problem is to adopt 1-irregular meshes and perform adaptive refinements toward the regions of interest. However, this strategy may lead to a high computational cost, as the number of degrees of freedom increases significantly. Alternatively, the directional refinement tool can be used to refine the mesh towards the stiffener line shown in Figure \ref{fig:casket-stiffners}. In Figure \ref{fig:casket-refined-mesh}, three steps of directional refinement are applied towards these lines, showing elements with good aspect ratio while creating graded elements in the direction of the region of interest.

\subsection{1-irregular and directional refinement of a cross-shaped geometry}

\begin{figure}[!ht]
  \centering
  \subfloat[]{\includegraphics[width=0.45\textwidth, valign=c]{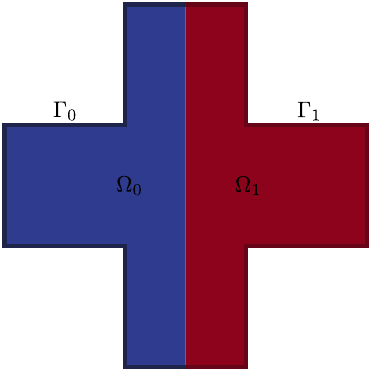}} \qquad
  \subfloat[\label{fig:cross-initial-mesh}]{\includegraphics[width=0.49\textwidth, valign=c]{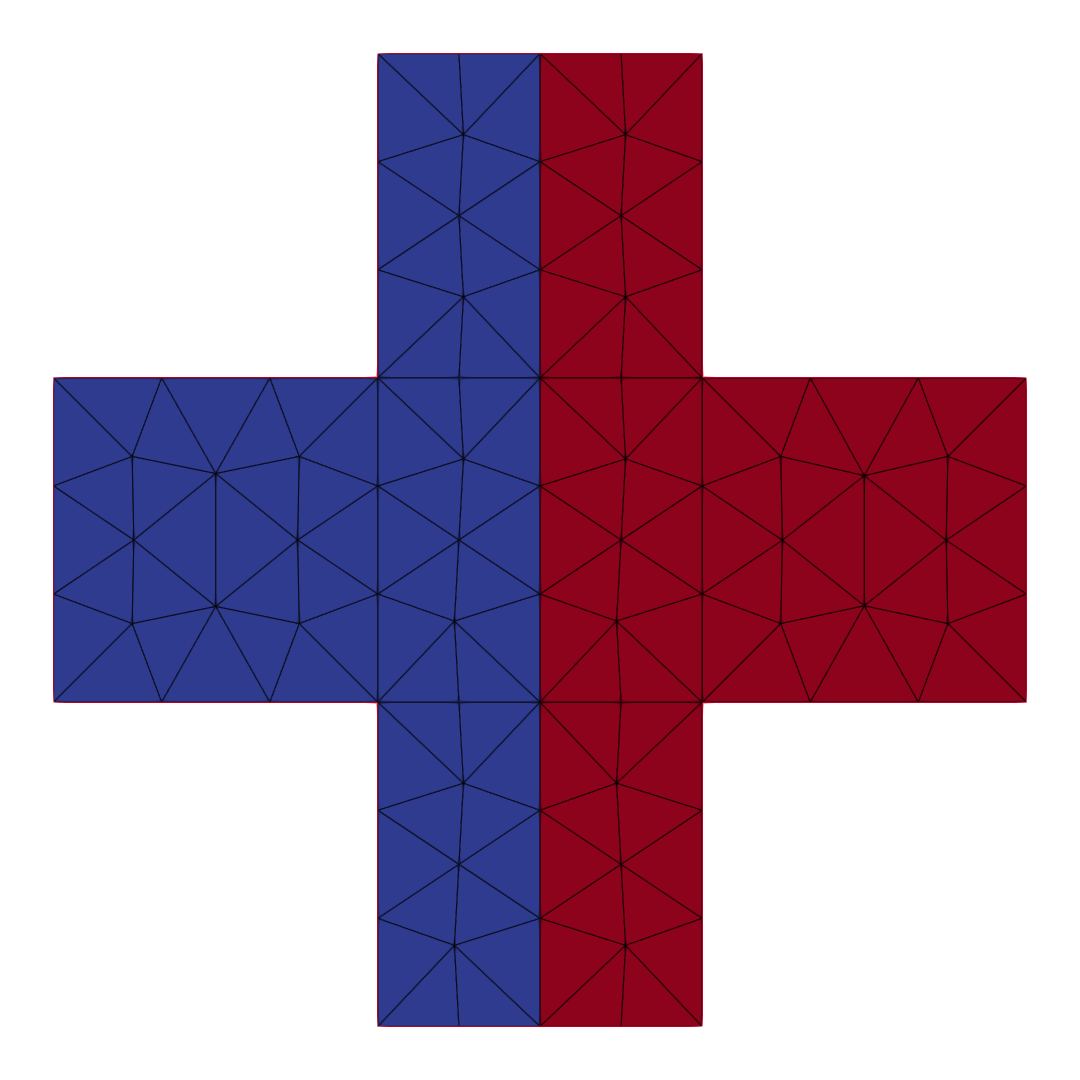}}
  \caption{Cross-shaped geometry. (a) subregions and boundary details and (b) initial triangular mesh}
  \label{fig:cross-initial-geometry}
\end{figure}

In fluid dynamics, the presence of a no-slip boundary condition (zero tangential velocity), often results in a \emph{boundary layer} scenario, in which the solution field presents a strong gradient towards the normal direction of the boundary of the domain. To analyze the applicability of the refinement patterns in such scenarios, Figure \ref{fig:cross-initial-geometry} illustrates a cross-shaped geometry divided into two subregions denoted by $\Omega_0$ and $\Omega_1$, whose boundaries are named $\Gamma_0$ and $\Gamma_1$, respectively.
For each subregion, a refinement strategy towards the boundary is applied. The elements of $\Omega_0$ are refined uniformly, \emph{i.e.}, using a refinement pattern in which the aspect ratio of the refined elements is similar to the original one. In $\Omega_1$, however, the \texttt{RefineDirectional} algorithm is applied. An unstructured mesh composed of triangular elements is adopted as an initial mesh for both subregions, as shown in Figure \ref{fig:cross-initial-mesh}.

\begin{figure}[!ht]
  \centering
  \subfloat[Level 0]{\includegraphics[width=0.49\textwidth]{cross-ref0-tri.png}} \hfill
  \subfloat[Level 1]{\includegraphics[width=0.49\textwidth]{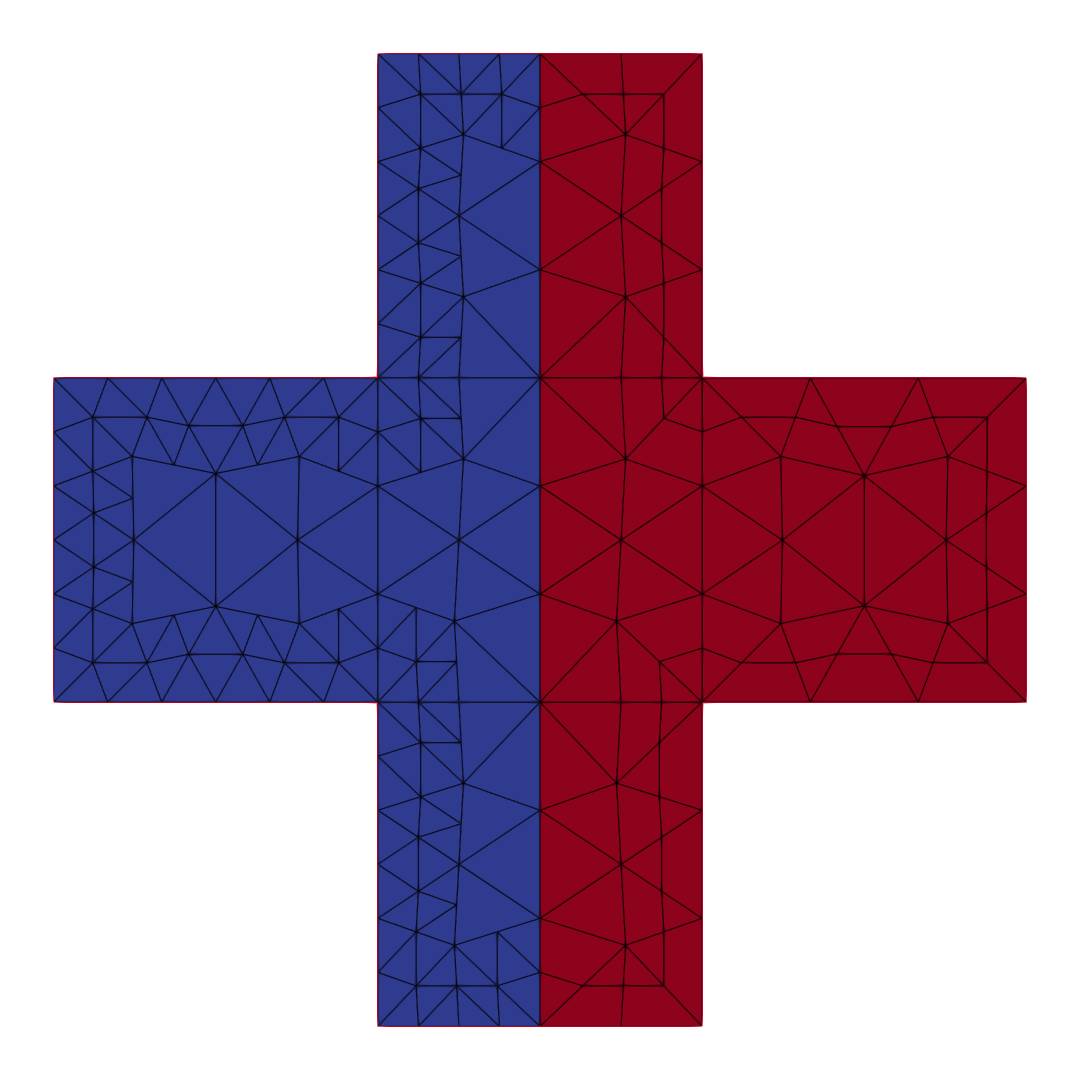}} \\
  \subfloat[Level 2]{\includegraphics[width=0.49\textwidth]{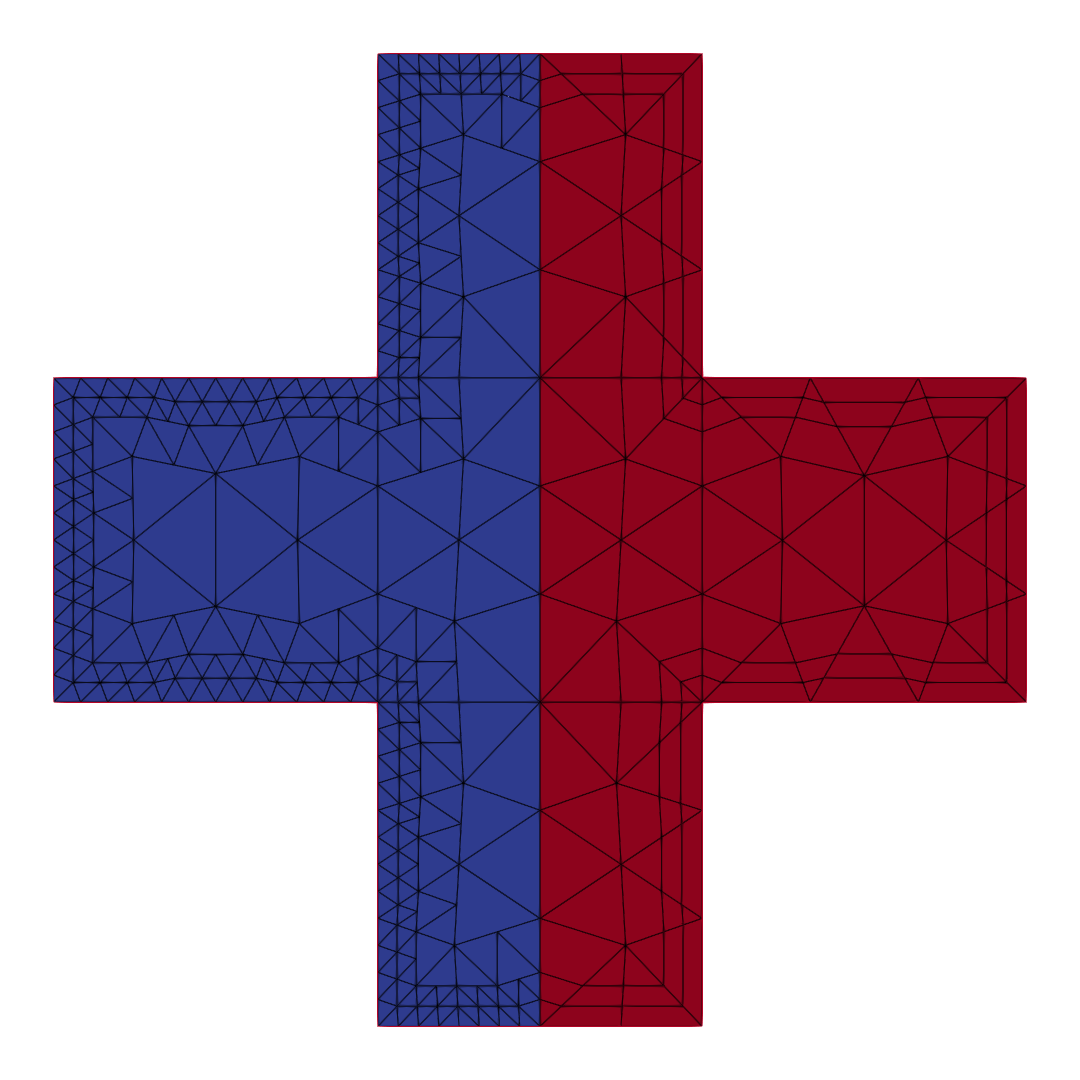}} \hfill
  \subfloat[Level 3]{\includegraphics[width=0.49\textwidth]{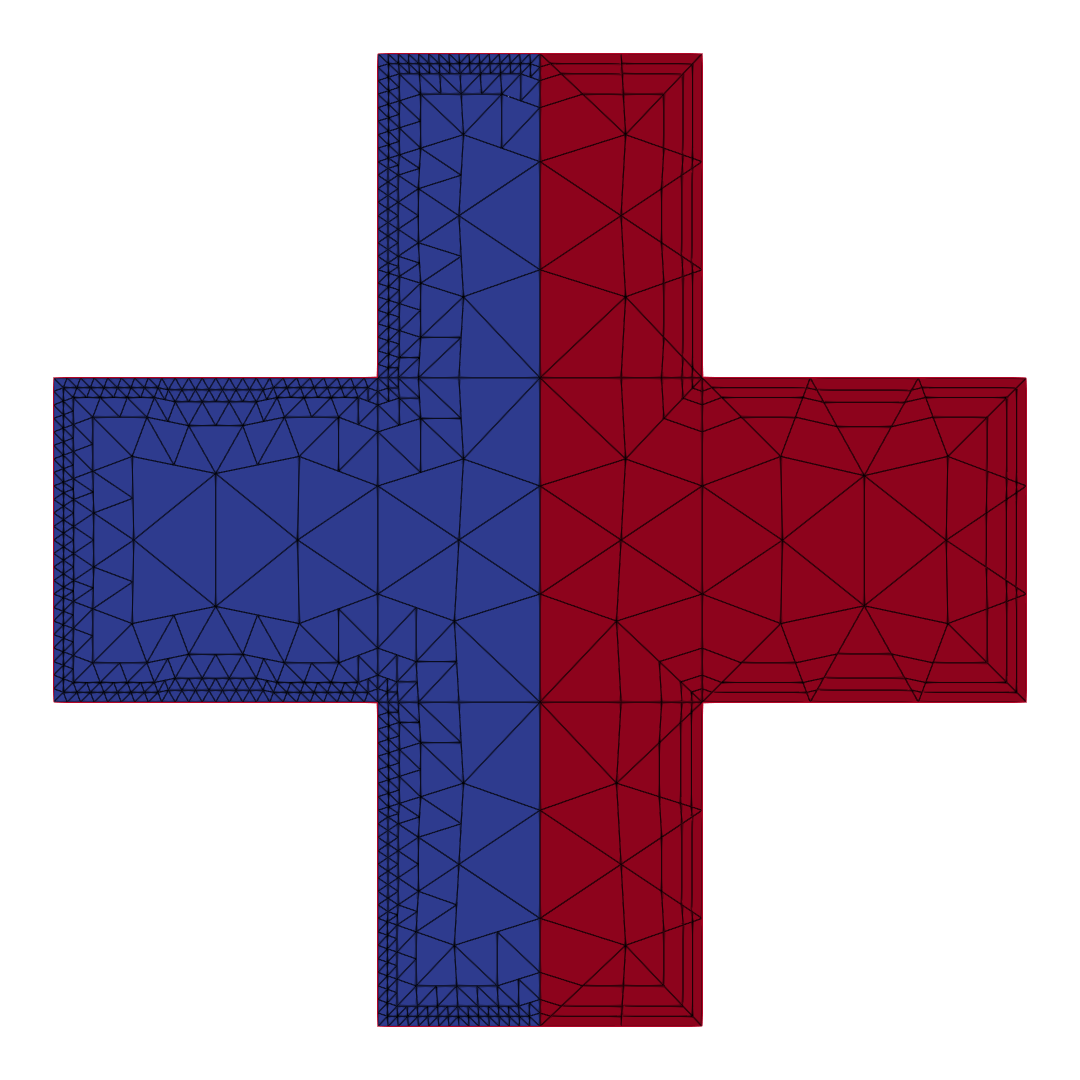}}
  \caption{Cross-shaped geometry. Sequence of refinements in both $\Omega_0$ and $\Omega_1$.}
  \label{fig:cross-refinements}
\end{figure}

The different refinement strategies are performed in each subregion in three subsequent levels. For each level, neighboring elements to the boundary are refined using the pattern adopted for that subregion. Figure \ref{fig:cross-refinements} shows the resulting mesh at each level of refinement.

The refinement strategy applied in $\Omega_1$ results in a significantly smaller number of elements, while retaining a similar element size in the normal direction to the boundary. Therefore, in a boundary-layer situation, this mesh is expected to be able to represent the gradients of the solution with better efficiency than with the uniform refinement strategy adopted in $\Omega_0$.

\section{Conclusions \label{sec:Conclusion}}

A methodology for generalizing h adaptivity in the context of Finite Element meshes is presented and the concept of refinement patterns is introduced. A refinement pattern is defined as an arbitrary partition of a master element. In order to describe refinement patterns formally, each master element is identified as the union of its sides. Each side of a sub-element is necessarily contained in a side of its father element, and an affine transformation allows to map points between sides. 

The refinement patterns documented in this work and their management in a database offer its users a toolset to refine a mesh in an innovative way. This is exemplified by the use of directional refinements which can be applied to a variety of problems where the solution field presents strong gradients in a specific direction. The effectiveness of the refinement patterns is demonstrated through specialized mesh refinements applied to problems with complex geometries. The methodology is implemented in the \texttt{NeoPZ} library, and the results obtained in the examples presented in this chapter show the potential of the proposed approach in improving the efficiency of finite element simulations.

\section{Extensions and future work \label{sec:FutureWork}}

The methodology presented in this work is a first step towards the development of a general framework for mesh adaptivity. Refinement patterns at this point assume fixed positions of the nodes in the master element. A natural extension of this work is to allow for the movement of nodes in the master element. This would allow for the creation of more general refinement patterns, leading to arbitrarily graded meshes.

It should also be noted that the methodology can be applied to higher-dimension meshes, generalizing the concept h-refinement.

\bigskip\noindent {\bf Acknowledgments:} The authors gratefully acknowledge the support of ANP - {\em Brazil's National Oil, Natural Gas and Biofuels Agency} through the R\&D levy regulation and TotalEnergies through research project 5955 with Unicamp. The authors also acknowledge the support of EPIC -- {\em Energy Production Innovation Center}, hosted by the University of Campinas (UNICAMP) and sponsored by Equinor Brazil and FAPESP -- {\em S\~ao Paulo Research Foundation} (2017/15736-3). Acknowledgments are extended to the Center for Petroleum Studies (CEPETRO) and the School of Civil Engineering, Architecture and Urban Planning (FECFAU). Author G. Avancini gratefully acknowledges Equinor/FAPESP and TotalEnergies/FUNCAMP for the financial support (grants 2023/06981-5 and 76042-23). Authors P.~R.~B. Devloo (grants 305823/2017-5 and  309597/2021-8) and F.~T. Orlandini (grant 130002/2022-7) thankfully acknowledges financial support from the CNPq - {\em Conselho Nacional de Desenvolvimento Cient\'\i fico e Tecnol\'ogico}.

\bibliographystyle{elsarticle-num} 
\addcontentsline{toc}{section}{\refname}
\bibliography{main}

\end{document}